\documentclass[11pt]{amsart}
\usepackage{preamble}

\title[Covariance of Error Terms]{Covariance of Error Terms Related to the Dirichlet Eigenvalue Problem }
\thanks{ This research was supported by the European Research Council (ERC) under the European Union's  Horizon 2020 research and innovation program  (Grant agreement No.    786758).\\
This work is part of the author’s M.Sc. thesis written under the supervision of Zeev Rudnick at Tel Aviv University.  \\
The author would like to extend his thanks to Igor Wigman and Nadav Yesha for their comments and corrections.}

\author{Noam Kimmel}

\begin{document}

\maketitle 
\begin{abstract}
    We explore the covariance of error terms coming from Weyl's conjecture regarding the number of Dirichlet eigenvalues up to size $X$.
    We also consider this problem in short intervals, i.e. the error term of the number of eigenvalues in the window $[X, X+S]$ for some $S(X)$.
    We look at these error terms for planar domains where the Dirichlet eigenvalues can be explicitly calculated.
    In these cases, the error term is closely related to the error term from the classical lattice points counting problem of expanding planar domains.
    We give a formula for the covariance of such error terms, for general planar domains.
    We also give a formula for the covariance of error terms in short intervals, for sufficiently large intervals.
    Going back to the Dirichlet eigenvalue problem, we give results regarding the covariance of the error terms in short intervals of 'generic' rectangles.
    We also explore a specific example, namely we compute the covariance between the error terms of an equilateral triangle and various rectangles.
\end{abstract}

\setcounter{tocdepth}{1}
\tableofcontents
\newpage

\section{Introduction}
\subsection{The Dirichlet Eigenvalue Problem}
\subsubsection{Definitions}
Let $\Gamma$ be a bounded planar domain with a piecewise smooth boundary $\partial \Gamma$.
The Dirichlet eigenvalues of $\Gamma$ are defined as the eigenvalues of (minus) the Laplace operator $\Delta$  corresponding to functions on $\Gamma$ which vanish on the boundary $\partial\Gamma$.
In other words, $\lambda$ is a Dirichlet eigenvalue if there is a non-zero solution $f$ to
\begin{gather*}
\begin{cases}
    -\Delta f = \lambda f \quad \text{in }\Gamma\\
    f|_{\partial\Gamma}=0 .         
\end{cases}    
\end{gather*}

It is known that there are infinitely many eigenvalues, and that these eigenvalues are positive and form a discrete set with a single accumulation point at infinity.
Thus, we can denote these eigenvalues as $0< \lambda_1\leq \lambda_2\leq ...\lambda_n\leq...$ where $\lim_{n\rightarrow\infty}\lambda_n = \infty$.

\subsubsection{Hearing the Shape of the Drum}
In 1966 Mark Kac published an article titled "Can One Hear the Shape of a Drum?" \cite{kac}.
This is a cute way of asking whether one can reconstruct the domain $\Gamma$ given its spectrum (since one can think of the original domain $\Gamma$ as a two dimensional membrane, and the eigenvalues are the possible frequencies of vibrations of that membrane where we hold the boundary in place).
In 1992 Carolyn Gordon, David Webb and Scott Wolpert found two different domains with the same spectrum \cite{G-W-W}, which showed that the general answer to the question is negative.
Nonetheless, it is still interesting to ask what information can be deduced about the original domain given its spectrum.

\subsubsection{Weyl's Law}
We denote $n(X) = \#\SET{\lambda_i\leq X}$ the cumulative counting function of the Dirichlet eigenvalues.
We also define $n(X,S) = n(X+S)-n(X)$ a window count function which
counts the number of Dirichlet eigenvalues in the segment $\left[X,X+S\right]$ where $S(X)$ depends on $X$.

An important theorem regarding these counting functions is Weyl's Law \cite{weyl}, which for the two dimensional case states that
$$n(X)\sim \frac{1}{4\pi}\text{Area}\left(\Gamma\right)X.$$
An immediate corollary from this theorem is that the area of the original domain $\Gamma$ can be deduced given its spectrum.
Weyl also conjectured the next term in the approximation-
$$
n(X) = \frac{1}{4\pi} \text{Area}(\Gamma) X - \frac{1}{4\pi} \text{Length}(\partial \Gamma) \sqrt{X} + e\left(X\right)
$$
where $e(x)$ is of smaller order.
While in general this conjecture is still open, it has been proven by Victor Ivrii \cite{ivrii} for domains $\Gamma$ where the set of periodic trajectories of a billiard in $\Gamma$ has measure 0.

\subsubsection{The Error Term and the Window Count Error Term}
There has been a lot of work dedicated to the study of $e(X)$ and other statistical properties of the spectrum.
For example, another important function is the error term for the window counting function, which can be expressed as $e(X,S)=e(X+S)-e(X)$.
The statistics of $e(X,S)$ behave differently depending on the growth rate of $S(X)$ as $X\rightarrow \infty$. 
The case where $S(X)$ is constant is especially important due to the Berry-Tabor conjecture (their original paper at \cite{berry-tabor}, and a survey can be found in \cite{zeev-QC}).
This is a major conjecture from the field of quantum chaos which states that if the corresponding classical dynamics of $\Gamma$ is completely integrable then the distribution $e(X,S)$ is usually Poisson.

\subsection{Lattice Point Counting Problem}
\subsubsection{Definitions}
We now consider another classical problem, that of estimating the number of lattice points inside an expanding planar domain.
Let $\Omega$ be a convex shape in the plane, where the boundary $\gamma = \partial \Omega$ is smooth and has positive curvature.
We denote 
$$
N(R) = \#\SET{\ZZ^2 \cap R\Omega}.
$$
It is known that $N(R)$ is aproximated by the area of $R\Omega$ which is $R^2\text{Area}\left(\Omega\right)$.
And so, we define the error function
$$E(R) = N(R) - \text{Area}\left(\Omega\right)R^2$$
\subsubsection{Bounding the Lattice Point Error Function}
There has been a lot of research dedicated to try and bound this error term $E(R)$.
One can easily show that $E(R) = \BigO{R}$.
For the case where $\Omega$ is the unit disc, this was shown by Gauss in 1834.
However, it was conjectured by Hardy that $E(R) = \BigO{R^{1/2 + \epsilon}}$ for all $\epsilon>0$.

While we still don't know to show that $E(R) = \BigO{R^{1/2 + \epsilon}}$, there have been significant improvements in this area.
The bound $E(R) = \BigO{R^{2/3}}$ was obtained by Sierpiński \cite{sierpinski} at the start of the twentieth century.
Since then, various improvements to this bound were obtained.
The current best known bound for the case of the circle is $E(R) = \BigO{R^{517/824+ \epsilon} }$ due to Bourgain and Watt \cite{B-W}.
The lower bound 
$\limsup_{R\rightarrow\infty}E(R)/R^{1/2} > 0$
was shown independently by Hardy and Landau in 1915 \cite{hardy} \cite{landau}.

\subsubsection{The Asymptotic Distribution of the Lattice Point Error Function}
While proving the bound $E(R) = \BigO{R^{1/2 + \epsilon}}$ is still outside our reach, in the 1940's Wintner considered a related problem \cite{wintner}.
Consider the normalized error function 
\begin{equation}\label{eq:normalized_error_function}
F(R) = E(R)/R^{1/2}.
\end{equation}
We can ask whether this function has an asymptotic density.
That is, we want to know if there exists a density function $p(t)$, such that for every segment $[a,b]\subset \RR$,
$$
\lim_{T\rightarrow\infty}\frac{1}{T}\mu\left(\SET{R\in[0,T] : F(R)\in[a,b]}\right)
=\int_{a}^{b}p(t)dt
$$
where $\mu$ is the standard Lebesgue measure.

Wintner and Heath-Brown showed that for the case of the circle such a density function does exist \cite{wintner} \cite{heath}. 
This was later generalized by Bleher for general convex domains with a boundary which is smooth and has positive curvature \cite{bleher}.

\subsection{The Connection Between Counting Dirichlet Eigenvalues and Counting Lattice Points}
\subsubsection{Domains with Explicit Eigenvalues Formulas}
There is a connection between the Dirichlet eigenvalue counting problem and the lattice point counting problem.
This connection becomes apparent once we consider the Dirichlet eigenvalue problem for domains where the eigenvalues can be explicitly calculated.
Very few such domains are known, and these include rectangles, ellipses, and a few select triangles such as the equilateral triangle.

Consider for example the case of a rectangle $\Gamma = [0,a]\times[0,b]$.
The Dirichlet eigenfunctions in this case are given by 
$$f(x,y)=\sin\left(\frac{\pi n x}{a}\right)
\sin\left(\frac{\pi m y}{b}\right), \quad n,m\in\NN$$
and the eigenvalues are
\begin{equation}\label{eq:rectangle_explicit_eigenvalues}
\SET{\pi^2\left(\frac{n^2}{a^2}+\frac{m^2}{b^2}\right) : n,m\in\NN}.    
\end{equation}
And so it can be seen that $n(X)$, the number of eigenvalues up to $X$, is related to the number of integer lattice points inside the ellipse 
$\Omega : \frac{x^2}{a^2} + \frac{y^2}{b^2}\leq \frac{1}{\pi^2}$
dilated by a factor of $\sqrt{X}$.
There are still minor differences between $N_{\Omega}(\sqrt{X})$ and $n_{\Gamma}(X)$ due to fact that the Dirichlet eigenvalues come from positive $n,m$.
Nonetheless, it can still be shown that
$$N_{\Omega}(X) = 4n_{\Gamma}\left(X^2\right) + \frac{2}{\pi}(a+b)X + \BigO{1}$$
and as a result 
\begin{equation}\label{eq:connection_ratio_rec}
E_{\Omega}(X) = 4e_{\Gamma}(X^2) +\BigO{1}.   
\end{equation}

A similar thing happens in the case of $\Gamma$ an equilateral triangle of side length $\ell$.
While not a trivial matter as in the case of the rectangle, it was shown by Lam\'e \cite{lame} that the eigenvalues of $\Gamma$ are of the form 
$$\frac{16\pi^2}{9 \ell^2}\left( n^2 + nm + m^2\right)$$
for $n,m\in \mathbb{N}$.
A nice account of this was given by McCartin \cite{McCartin}.
Denoting by $\Omega$ the ellipse $x^2 + xy + y^2 \leq \frac{9 \ell^2}{16\pi^2}$ we once more get a connection
$$N_{\Omega}(X) = 6n_{\Gamma}(X^2) + \frac{9}{2\pi}\ell \sqrt{X} + \BigO{1}.$$
which implies 
\begin{equation}\label{eq:connection_ratio_tri}
E_{\Omega}(X) = 6e_{\Gamma}(X^2) + \BigO{1}.
\end{equation}

We see that for these cases, the study of the error term from the Dirichlet eigenvalue problem is equivalent to the study of the error term from the lattice point problem for ellipses.

\subsubsection{The Window Count of Dirichlet Eigenvalues and the Lattice Point Problem in Thin Annuli}
In the case of rectangles and special triangles, there is also an analogue for the window count function in terms of the lattice point problems.
The number of lattice points in a thin annuli of shape $\Omega$ and width $h$, is given by $N_{\Omega}(R,h)=N_{\Omega}(R+h)-N_{\Omega}(R)$.
The main term of $N_{\Omega}(R,h)$ as $R$ goes to infinity is $\text{Area}\left(\Omega\right)(2hR+h^2)$.
The error term is given by $E(R,h)=E(R+h)-E(R)$, and we define the normalized error term as $F(R,h)=E(R,h)/R^{1/2}$.
Note that in the case of $e(X,S)$, $S$ is the area of the annulus,
while in $E(R,h)$, $h$ is its width.

\subsection{Main Results}
In this paper we investigate further statistical properties of the error functions $E(X)$ and $E(X,h)$ coming from the lattice point problem for convex smooth domains.
Specifically, we consider the covariance of such error function coming from two different domains.
We give a more precise formula for the case of ellipses.
These cases are of particular importance as they are also equivalent to error functions coming from the Dirichlet eigenvalue problem.

\subsubsection{Global Covariance}
Let $\Omega_1,\Omega_2$ be two convex planar domains, each with a boundary which is smooth and has positive curvature, and let $F_1,F_2$ be their normalized error functions for the lattice counting problem as defined in \eqref{eq:normalized_error_function}.
It was shown by Bleher that these functions have an asymptotic density, and he also gave a formula for their variance.
We explore the covariance of these functions, where the covariance of $F_1$ and $F_2$ is given by the following definition.
\begin{definition}[Covariance] \label{defn:covar}
$$
\COV{F_1,F_2} = \lim_{T\rightarrow\infty}
\frac{1}{T}\int_{0}^{T}F_1(t)F_2(t)dt .
$$
\end{definition}

In \autoref{thm:general_global_covar} we give a formula for $\COV{F_1,F_2}$.
\begin{theorem} \label{thm:general_global_covar}
Let $\Omega_1,\Omega_2$ be two convex planar domains, each with a boundary which is smooth and has positive curvature.
Denote $\gamma_i=\partial \Omega_i$, $i=1,2$.
For $v\in\RR^2$ we denote $x_i(v)$ the point on $\gamma_i$ with outer normal $\frac{v}{\|v\|}$.
We define $Y_i(v) = v\cdot x_i(v)$ and $\rho_i(v)$ as the curvature radius of $\gamma_i$ at $x_i(v)$.
With these notations-
$$
\COV{F_1,F_2} = \frac{1}{2\pi^2}\sideset{}{'}\sum_{n\in\ZZ^2}\sideset{}{'}\sum_{\substack{m\in\ZZ^2\\ Y_2(m) = Y_1(n)}}
\frac{\sqrt{\rho_1(n)}\sqrt{\rho_2(m)}}{|n|^{3/2} |m|^{3/2}}
$$
where $\sideset{}{'}\sum$ denotes a sum where $0$ is omitted.
\end{theorem}
\begin{remark}
A typical pair of domains $\Omega_1,\Omega_2$ will have no non-zero $n,m\in\ZZ^2$ such that $Y_1(n) = Y_2(m)$ (which means that the covariance will be 0).
If we restrict to the case where the covariance is positive, then for the typical case the set of solutions $n,m\in\ZZ^2$ to  $Y_1(n) = Y_2(m)$ will be of the form $\SET{(k n_0, k m_0) | k\in \ZZ}$ for some initial solution $n_0,m_0\in\ZZ^2$.

However, regardless of the number of summands, the series in \autoref{thm:general_global_covar} always converges.
This is a result of the function $F_i$ being in the space $\mathbb{B}^2$ of almost periodic functions, as will be discussed in \autoref{sec:two}.
\end{remark}

Applying \autoref{thm:general_global_covar} for the case of ellipses, we get the following corollary:
\begin{corollary}\label{coro:ellipse_global_covar}
Let $\Omega_1$, $\Omega_2$ be ellipses of the form 
$$
\Omega_1: ax^2 + bxy + cy^2 \leq 1,\quad \Omega_2: dx^2 + exy + fy^2 \leq 1
$$
For some parameters $a,b,c,d,e,f$.
Then the covariance between the normalized error functions $F_{\Omega_1}(t), F_{\Omega_2}(t)$ is 
$$
\frac{16\pi}{\sqrt{4ac-b^2}\sqrt{4df-e^2}}
\sum_{\substack{\nu:\\ r_1(\nu)\neq 0, r_2(\nu)\neq 0}}
\frac{r_1(\nu)r_2(\nu)}{\nu^{3}}
$$
where $r_1(\nu)$ is the number of integer solutions $n_1,n_2$ to 
$$
4\pi
\frac{\sqrt{c n_1^2 - b n_1 n_2 + an_2^2}}{\sqrt{4ac-b^2}} = \nu
$$
and similarly $r_2(\nu)$ is the number of integer solutions $m_1,m_2$ to 
$$
4\pi
\frac{\sqrt{f m_1^2 - e m_1 m_2 + d m_2^2}}{\sqrt{4df-e^2}} = \nu.
$$
\end{corollary}
\begin{remark}
Note that the sets $\left\lbrace \nu : r_i(\nu)\neq 0\right\rbrace$, $i=1,2$, are discrete. 
Thus, the sum in \autoref{coro:ellipse_global_covar} over $\nu$ such that $r_1(\nu)\neq 0, r_2(\nu)\neq 0$ is a sum over a discrete set.
\end{remark}

\subsubsection{Covariance in Short Intervals}
We also consider the covariance of the error term for short intervals.
Let $F_1(R,h),F_2(R,h)$ be the normalized error functions in short intervals of the lattice point problem related to $\Omega_1,\Omega_2$.
We assume that $h(R)$ tends to 0 as $R$ tends to infinity.
We define the covariance of these functions as
\begin{definition}[Covariance in Short Intervals]\label{defn:covariance_in_short_intervals}
For a given function $h(R)$, we say that the covariance between $F_1(R,h)$ and $F_2(R,h)$ is $f(h)$ for some function $f$ if 
$$
\frac{1}{T}\int_{t=0}^{T}F_1(t,h(T))F_2(t,h(T))dt \sim f(h(T)).
$$
\end{definition}
In \autoref{thm:covar_short_interval_F1_F2} we give a formula for the covariance in short intervals.
However, our proof techniques require us to have $h(R)$ decay sufficiently slowly towards 0.
Furthermore, we also get a Diophantine-like condition that the two domains need to satisfy.
For general domains this condition is somewhat vague. 

We can however make these conditions more concrete if we restrict our attention to a small space of domains.
In \autoref{chap:ellipse_short} we restrict our attention to ellipses.
In this case, Bleher and Lebowitz showed \cite{bleher_lebowitz} that the variance of the normalized error function $F(R,h)$ of a "generic" ellipse of the form $\SET{(x,y) : x^2 + \mu y^2 \leq 1}$ has order $h$, provided that $h(T) \gg T^{-1 + \epsilon}$ for some $\epsilon>0$.
We show that for a generic pair of ellipses that have non-zero covariance in the global case, the covariance in a short interval $h$ will be of order $h^2\log\left( h^{-1}\right)$ provided that $h(T) \gg T^{-1/110 + \epsilon}$ for some $\epsilon>0$.
The main result for this is
\begin{theorem}\label{thm:ellipse_short_covar}
Almost all pairs of ellipses
$$
\Omega_1: ax^2 + bxy + cy^2 \leq 1,\quad
\Omega_2: dx^2 + exy + fy^2 \leq 1
$$
that have non-zero global covariance satisfy the following:

For any $h(T)$ with $h(T) \gg T^{-1/110 + \epsilon}$ for some $\epsilon>0$:
$$
\COV{F_1(T,h),F_2(T,h)} = C h^2\log\left( h^{-1}\right)
$$
where $F_1,F_2$ are the normalized error functions of $\Omega_1, \Omega_2$.

The constant $C$ is given by 
$$
C = \frac{16\pi}{\nu_0\sqrt{4ac-b^2}\sqrt{4df-e^2}}
$$
where $\nu_0$ is the smallest positive number for which
$$
4\pi
\frac{\sqrt{c n_1^2 - b n_1 n_2 + an_2^2}}{\sqrt{4ac-b^2}}
= 
4\pi
\frac{\sqrt{f m_1^2 - e m_1 m_2 + d m_2^2}}{\sqrt{4df-e^2}}
= \nu_0
$$
has an integer solution $n_1,n_2,m_1,m_2$.
\end{theorem}

\begin{remark}
The precise definition of "almost all" from \autoref{thm:ellipse_short_covar} will be given in \autoref{sec:gen_elips_short_intervals}.
\end{remark}

We also examine the covariance between two ellipses coming from the Dirichlet eigenvalue problem of an equilateral triangle and a rectangle, for an interesting family of rectangles.
We get the following theorem:

\begin{theorem}\label{thm:specific_ellipse}
Let $F_{1},F_{2}$ be the normalized error functions of the following two ellipses:
$$\Omega_1: x^2 + xy + y^2 \leq \frac{3}{4}, \quad \Omega_2 : x^2 + \alpha y^2 \leq 1$$
for some parameter $\alpha\in\RR^+$.
\begin{enumerate}
\item For almost all $\alpha\in\RR^+$, if $h(T)\gg T^{-1/86+\epsilon}$ for some $\epsilon>0$, then the covariance of $F_{1}$ and $F_{2}$ is $\frac{9}{\pi\sqrt{\alpha}} h^2 \log^2 \left(h^{-1}\right)$

\item If $\alpha=\frac{p}{q}$ is a rational such that $3pq$ is not a perfect square, then assuming $h(T) \gg T^{-1/62}$, the covariance of $F_{1}$ and $F_{2}$ is $\frac{C \sqrt{3}}{\sqrt{\alpha}}h$ where $C$ is a constant given in \autoref{thm:app_non_square_case}.

\item If $\alpha=\frac{p}{q}$ is a rational and $3pq$ is a perfect square then assuming $h(T) \gg T^{-1/62}$, the covariance of $F_{1}$ and $F_{2}$ is $\frac{18}{\sqrt{pp'}}h\log\left(h^{-1}\right)$
where $p'$ is the squarefree part of $p$.
\end{enumerate}
\end{theorem}

Using \eqref{eq:connection_ratio_rec},\eqref{eq:connection_ratio_tri},
we can restate 
\autoref{thm:specific_ellipse} in terms of the Dirichlet eigenvalue problem:

\begin{corollary}
Let $\Gamma_1$ be an equilateral triangle of side length $\frac{2\pi}{\sqrt{3}}$, and let $\Gamma_2$ be a rectangle with side lengths $\pi$ and $\sqrt{\alpha}\pi$ for some $\alpha>0$.
Denote by $n_i(t)$ the number of Dirichlet eigenvalues of $\Gamma_i$ up to $t^2$, and denote
$$
e_i(t) = n_i(t) - \frac{1}{4\pi}\text{Area}(\Gamma_i)t^2
+  \frac{1}{4\pi}\text{Length}(\partial\Gamma_i) t.
$$
\begin{enumerate}
\item For almost all $\alpha$, if  $h(X)\gg X^{-1/86+\epsilon}$ for some $\epsilon>0$, then
$$
\frac{1}{X}\int_{0}^{X}
\frac{\left(e_1(t+h)-e_1(t)\right)\left(e_2(t+h)-e_2(t)\right)}{t}dt
\sim \frac{3}{8\pi\sqrt{\alpha}} h^2 \log^2 \left(h^{-1}\right).
$$
\item If $\alpha=\frac{p}{q}$ is a rational such that $3pq$ is not a perfect square, and if $h(T) \gg X^{-1/62}$, then 
$$
\frac{1}{X}\int_{0}^{X}
\frac{\left(e_1(t+h)-e_1(t)\right)\left(e_2(t+h)-e_2(t)\right)}{t}dt
\sim \frac{C }{8\sqrt{3\alpha}}  h
$$
where the constant $C$ is given in \autoref{thm:app_non_square_case}.
\item If $\alpha=\frac{p}{q}$ is a rational and $3pq$ is a perfect square, and if $h(T) \gg X^{-1/62}$, then
$$
\frac{1}{X}\int_{0}^{X}
\frac{\left(e_1(t+h)-e_1(t)\right)\left(e_2(t+h)-e_2(t)\right)}{t}dt
\sim \frac{3}{4 \sqrt{pp'}}h\log\left(h^{-1}\right)
$$
where $p'$ is the squarefree part of $q$.
\end{enumerate}
\end{corollary}

Using \eqref{eq:rectangle_explicit_eigenvalues}, \autoref{thm:ellipse_short_covar} can also be restated in terms of the error terms of the Dirichlet eigenvalue problem.
This is not a direct corollary, but similar arguments to those given in \autoref{sec:gen_elips_short_intervals} result in the following:

\begin{corollary}
Let $\Gamma_1 = [0,a]\times [0,b],\Gamma_2=[0,c]\times[0,d]$ be two rectangles.
Denote by $n_i(t)$ the number of Dirichlet eigenvalues of $\Gamma_i$ up to $t^2$, and denote
$$
e_i(t) = n_i(t) - \frac{1}{4\pi}\text{Area}(\Gamma_i)t^2
+  \frac{1}{4\pi}\text{Length}(\partial\Gamma_i) t.
$$
For almost all pairs of rectangles $\Gamma_1,\Gamma_2$ that have non-zero global covariance, we have that if $h(T) \gg T^{-1/110 + \epsilon}$ for some $\epsilon>0$:
$$
\frac{1}{X}\int_{0}^{X}\frac{\left(e_1(t+h)-e_1(t)\right)\left(e_2(t+h)-e_2(t)\right)}{t}dt
\sim
C h^2\log\left( h^{-1}\right).
$$
The constant $C$ is given by 
$$
C = \frac{a b c d}{(2\pi)^3 \nu_0}
$$
where $\nu_0$ is the smallest positive number for which
$$
 \sqrt{b^2 n_1^2 + a^2 n_2^2}
= 
 \sqrt{d^2 m_1^2 + c^2 m_2^2}
= \nu_0
$$
has an integer solution $n_1,n_2,m_1,m_2$.
\end{corollary}

Similar results can also be obtained for a rectangle and equilateral triangle pair.

The reason we restrict our attention to rectangles is that their Dirichlet eigenvalues can be explicitly calculated.
For general parallelograms the eigenvalues with Dirichlet or Neumann boundary conditions are not explicitly known (only with periodic boundary conditions).
In fact, for parallelograms which are not rectangles, the Dirichlet eigenfunctions are not trigonometric, see \cite[Theorems 4.1.2, 4.2.1]{McCartin}.

\section{Global Covariance}\label{sec:two}

Let $\Omega_1,\Omega_2$ be two convex planar domains, each with a boundary which is smooth and has positive curvature.
In this section we look at the global normalized error functions $F_1(R),F_2(R)$ of these domains as defined in \eqref{eq:normalized_error_function}.
In order to find the covariance of these functions we use the fact that these functions are Besicovitch almost periodic functions.
In \autoref{sec:almost_periodic_functions} we introduce the theory of almost periodic functions.
In \autoref{sec:global_covariance} we apply the theory to the normalized error functions to obtain a formula for the global covariance.
In \autoref{sec:ellipse_golbal_covariance} we consider the more specific case where $\Omega_1,\Omega_2$ are ellipses.

\subsection{Almost Periodic Functions}\label{sec:almost_periodic_functions}

\subsubsection{The space $\BB^2$}
\begin{definition}\label{def:bb2}
The space $\BB^2$ of Besicovitch almost periodic functions (see \cite{besicovitch}) is defined as the closure of trigonometric polynomials under the semi-norm
$$
\|f\| = \left(\limsup_{T\rightarrow \infty}\frac{1}{T}\int_{0}^{T}\left|f\right|^2(t)dt\right)^{1/2}.
$$

In other words, a function $f:\RR\rightarrow\RR$ is in $\BB^2$ if there exists a sequence of trigonometric polynomials $f_N(t) = \sum_{n\leq N}c_n e^{i\lambda_n t}, (c_n\in\mathbb{C},\lambda_n\in\RR)$ such that 
$$
\lim_{N\rightarrow\infty}\limsup_{T\rightarrow \infty}\frac{1}{T}\int_{0}^{T}\left|f-f_N\right|^2(t)dt = 0.
$$
\end{definition}

For a function $f\in\BB^2$, the limit
$$
\lim_{T\rightarrow\infty}\frac{1}{T}\int_{0}^{T}\left|f\right|^2(t)dt
$$
exists.

\subsubsection{Fourier Transform of Almost Periodic Functions}
\begin{definition}
The Fourier transform of a function $f\in\BB^2$ is defined as
$$
\hat{f}(\xi) = 
\lim_{T\rightarrow \infty}\frac{1}{T}\int_{0}^{T}f(t)e^{- i \xi t}dt .
$$
For functions $f\in\BB^2$, this limit exists for all $\xi\in\RR$.
\end{definition}

An important result in this regard is the following lemma.
\begin{lemma}
$\hat{f}(\xi)$ is 0 for all $\xi$ except at most a countable set.
\end{lemma}

If $\hat{f}(\xi_0)\neq 0$, we refer to $\xi_0$ as a frequency of $f$.

\begin{definition}[Frequency Expansion]
If $\xi_1,\xi_2, ...$ are all the frequencies of $f\in\BB^2$, we denote 
$$\sum_{n=1}^{\infty}\hat{f}(\xi_n) e^{i\xi_n t}$$
as the frequency expansion of $f$.
\end{definition}

Note that this frequency expansion can be a divergent series.
As such, it should be regarded as a formal way of encoding the Fourier transform of $f$, and not as a formula.

\subsubsection{Parseval's Identity for Almost Periodic Functions}
Our main gain from introducing Fourier transform of functions in $\BB^2$ is the following theorem.
\begin{theorem}[Parseval's Identity]\label{thm:parseval}
Let $f$ be a function in $\BB^2$, and let $\sum_{n=1}^{\infty}\hat{f}(\xi_n) e^{i\xi_n t}$ be its frequency expansion.
Then
$$
\|f\|^2 = \lim_{T\rightarrow \infty}\frac{1}{T}\int_{0}^{T}\left|f\right|^2(t)dt
=
\sum_{n=1}^{\infty}\left|\hat{f}(\xi_n)\right|^2.
$$
\end{theorem}

\begin{corollary}\label{coro:parseval}
Let $f,g\in\BB^2$ with frequency expansions
$$
\sum_{n=1}^{\infty}\hat{f}(\alpha_n) e^{i\alpha_n t},\quad 
\sum_{n=1}^{\infty}\hat{g}(\beta_n) e^{i\beta_n t}.
$$
Denote $\SET{\xi_n} = \SET{\alpha_n}\cap\SET{\beta_n}$ the set of common frequencies.
Then
$$
\lim_{T\rightarrow\infty}\frac{1}{T}
\int_{0}^{T}f(t)\overline{g(t)}dt =
\text{Re}\left(\sum_{n=1}^{\infty}\hat{f}(\xi_n) \overline{\hat{g}(\xi_n)}\right).
$$
\end{corollary}
\begin{proof}
Apply \nameref{thm:parseval} to $f+g$, which is also in $\BB^2$.
Note that the frequency expansion of $f+g$ is the sum of the frequency expansions of $f$ and $g$.
\end{proof}

\subsection{Covariance of General Planar Domains}\label{sec:global_covariance}
\subsubsection{Definitions}
We reintroduce some of the definitions given in \autoref{thm:general_global_covar}.
Let $\Omega$ be a convex planar domain with a boundary $\gamma = \partial\Omega$ which is smooth and has positive curvature.
For $v\in\RR^2$ we denote $x(v)$ the point on $\gamma$ with outer normal $\frac{v}{\|v\|}$.
We define $Y(v) = v\cdot x(v)$ and $\rho(v)$ as the curvature radius of $\gamma$ at $x(v)$.
Denote also $F(R)$ the normalized error function associated with $\Omega$ as defined in \eqref{eq:normalized_error_function}.

\subsubsection{Almost Periodicity of the Error Function}
In \cite{bleher} Bleher proved that the normalized error function $F(R)$ is in $\BB^2$.
More specifically, we have
\begin{theorem} \label{thm:trig_aprox}
Denote
\begin{multline*}
P_N(t) = \\
\frac{1}{2\pi}\sum_{\substack{n\in\ZZ^2 \\ 0<|n|<N}}\frac{\sqrt{\rho(n)}}{|n|^{3/2}}
\left(
\exp \left(i2\pi Y(n)t-i\frac{3\pi}{4} \right)
+ 
\exp \left(-i2\pi Y(n)t+i\frac{3\pi}{4} \right)
\right)
\end{multline*}
which can be more compactly written as
$$
P_N(t) = 
\frac{1}{\pi}\sum_{\substack{n\in\ZZ^2 \\ 0<|n|<N}}\frac{\sqrt{\rho(n)}}{|n|^{3/2}}
\cos\left(2\pi Y(n) t - \frac{3\pi}{4}\right).
$$
For $N\in \NN,T\in \RR$ such that $N^{5/2}\leq T$ we have
$$
\frac{1}{T}\int_{0}^{T} \left|F(t)-P_N(t)\right|^2dt
 \ll_N N^{-1/3}.
$$
\end{theorem}
\begin{corollary}
$$
\lim_{N\rightarrow\infty}\limsup_{T\rightarrow \infty}\frac{1}{T}\int_{0}^{T}\left|F(t)-P_N(t)\right|^2dt = 0.
$$
which implies $F(R)\in\BB^2$.
\end{corollary}

\begin{proof}[proof of \autoref{thm:trig_aprox}]
see \cite{bleher}, Theorem 3.1 and Lemmas 3.2, 3.3, 3.4.
These Lemmas state that there exists constants $C_1,C_2,C_3$ such that
$$
\frac{1}{T}\int_{0}^{T} \left|F(t)-P_N(t)\right|^2dt
\leq C_1 T^{-1/3} + C_2 \left(N^{-1/3} + \frac{\log^2 T}{T}\right) + C_3 T^{-2}.
$$
Since $N^{5/2}\leq T$, this completes the proof.
\end{proof}

Furthermore, it can be shown that the frequency expansion of $F(R)$ is
\begin{equation}\label{eq:frequency_expansion}
\frac{1}{2\pi}
\sideset{}{'}\sum_{n\in\ZZ^2 }
\frac{\sqrt{\rho(n)}}{|n|^{3/2}}
\left(
\exp \left(i2\pi Y(n)t-i\frac{3\pi}{4} \right)
+ 
\exp \left(-i2\pi Y(n)t+i\frac{3\pi}{4} \right)
\right).
\end{equation}

\subsubsection{A Formula For the Covariance}
Consider once more the problem of computing the covariance of $F_1(R),F_2(R)$, the normalized error functions associated with $\Omega_1$ and $\Omega_2$.

\begin{proof}[Proof of \autoref{thm:general_global_covar}]
From \eqref{eq:frequency_expansion} we have an expression for the frequency expansions of $F_1$ and $F_2$.
Furthermore, from \autoref{coro:parseval} we get a formula for the covariance of two functions in $\BB^2$.
Putting these together we get that
$$
\COV{F_1,F_2} = \frac{1}{2\pi^2}\sideset{}{'}\sum_{n\in\ZZ^2}\sideset{}{'}\sum_{\substack{m\in\ZZ^2\\ Y_2(m) = Y_1(n)}}
\frac{\sqrt{\rho_1(n)}\sqrt{\rho_2(m)}}{|n|^{3/2} |m|^{3/2}}
$$
which proves \autoref{thm:general_global_covar}.
\end{proof}

\begin{remark}
Note that the covariance is always non-negative, and it is positive if and only if $F_1$ and $F_2$ share common frequencies.
\end{remark}

\begin{corollary}
Let $\Omega_1,\Omega_2$ be two convex planar domains, each with a boundary which is smooth and has positive curvature.
The covariance between the normalized error functions of $\Omega_1$ and $\alpha\Omega_2$ is 0 except for a countable set of $\alpha's$ for which the covariance is positive.
\end{corollary}
\begin{proof}
The covariance is positive if and only if the normalized error functions share common frequencies.
This happens only when $\alpha$ is of the form 
$$
\alpha = Y_1(n) / Y_2(m), \qquad 0\neq n,m\in \ZZ^2.
$$
\end{proof}

\subsection{Covariance of Ellipses}\label{sec:ellipse_golbal_covariance}
\subsubsection{Calculations Related to Ellipses}
We can get a more specific formula than that in \autoref{thm:general_global_covar} if we restrict to the case of ellipses.
Let $\Omega$ be a general ellipse given by
$$
\Omega : ax^2+bxy+cy^2\leq 1.
$$

Denote $H(x,y) = ax^2+bxy+cy^2-1$.
The normal at a point $(x_0,y_0)$ on $\gamma=\partial\Omega$ is given by $$\left(H_x(x_0,y_0),H_y(x_0,y_0)\right) = \left(2ax_0+by_0,2cy_0+bx_0\right).$$
It follows that
\begin{multline*}
x(n_1,n_2) = \\
\left(\frac{2cn_1-bn_2}{\sqrt{4ac-b^2}\sqrt{an_2^2-bn_1n_2+cn_1^2}},\frac{2an_2-bn_1}{\sqrt{4ac-b^2}\sqrt{an_2^2-bn_1n_2+cn_1^2}}\right)
\end{multline*}

and
\begin{equation}\label{eq:y_ellipse}
Y(n_1,n_2) = (n_1,n_2)\cdot x(n_1,n_2) = 
\frac{2\sqrt{an_2^2-bn_1n_2+cn_1^2}}{\sqrt{4ac-b^2}}.
\end{equation}

The formula for the radius of curvature at a point $(x_0,y_0)$ on $\gamma$ is
$$
\frac{(H_x^2 + H_y^2)^{3/2}}{H_{xx}H_y^2 - 2H_{xy}H_x H_y + H_{yy}H_x^2}(x_0,y_0).
$$
It follows that
\begin{equation}\label{eq:rho_ellipse}
\rho(n_1,n_2) = \frac{\sqrt{4ac-b^2}}{2}\left(\frac{n_1^2+n_2^2}{an_2^2-bn_1n_2+cn_1^2}\right)^{3/2}.
\end{equation}

It will be useful to note that
\begin{multline}\label{eq:ellipse_freq_coef}
\frac{1}{2\pi}
\frac{\sqrt{\rho(n_1,n_2)}}{\left(n_1^2 + n_2^2\right)^{3/4}}
=
\frac{\sqrt[4]{4ac-b^2}}{2\sqrt{2}\pi\left(an_2^2-bn_1n_2+cn_1^2\right)^{3/4}}
=\\
\frac{2 \sqrt{2\pi}}{\sqrt{4ac-b^2}}
\left(
\frac{4ac-b^2}{16\pi^2\left(an_2^2-bn_1n_2+cn_1^2\right)}
\right)^{3/4}
=  \\
\frac{\sqrt{8\pi}}{\sqrt{4ac-b^2}}\left(2\pi Y(n_1,n_2)\right)^{-3/2}.
\end{multline}

\subsubsection{A Formula For The Covariance of Ellipses}

We can now prove the formula for the covariance of normalized error functions of ellipses.
\begin{proof}[Proof of \autoref{coro:ellipse_global_covar}]

Let $\Omega_1$, $\Omega_2$ be general ellipses of the form 
$$
\Omega_1: ax^2 + bxy + cy^2 \leq 1,\quad \Omega_2: dx^2 + exy + fy^2 \leq 1
$$
for some parameters $a,b,c,d,e,f$.
We plug \eqref{eq:y_ellipse} and \eqref{eq:rho_ellipse} into the formula of the covariance from \autoref{thm:general_global_covar}.
We get that
\begin{multline*}
\COV{F_1,F_2} = 
\frac{1}{2\pi^2}
\sideset{}{'}\sum_{n\in\ZZ^2}
\sideset{}{'}\sum_{\substack{m\in\ZZ^2\\ Y_2(m) = Y_1(n)}} \frac{\sqrt{\rho_1(n)}\sqrt{\rho_2(m)}}{|n|^{3/2} |m|^{3/2}} 
=  \\
2
\sideset{}{'}\sum_{n\in\ZZ^2}
\sideset{}{'}\sum_{\substack{m\in\ZZ^2\\ Y_2(m) = Y_1(n)}}
\left(
\frac{\sqrt{\rho_1(n)}}{2\pi |n|^{3/2}}
\right)
\left(
\frac{\sqrt{\rho_2(m)}}{2\pi |m|^{3/2}}
\right).
\end{multline*}

From \eqref{eq:ellipse_freq_coef}, we get that this is equal to
\begin{equation*}
\COV{F_1,F_2} = 
\sideset{}{'}\sum_{n\in\ZZ^2}
\sideset{}{'}\sum_{\substack{m\in\ZZ^2\\ Y_2(m) = Y_1(n)}}
\frac{16\pi}{\sqrt{4ac-b^2}\sqrt{4df-e^2}}
\frac{1}{(2\pi Y_1(n))^{3}}.
\end{equation*}

Denoting by $r_1(\nu)$ the number of integer solutions $n_1,n_2$ to
$$
4\pi \frac{\sqrt{a n_2^2 -bn_1n_2 +cn_1^2 }}{\sqrt{4ac-b^2}} = \nu
$$
and similarly $r_2(\nu)$ the number of integer solutions $m_1,m_2$ to
$$
4\pi \frac{\sqrt{d m_2^2 -e m_1 m_2 +f m_1^2 }}{\sqrt{4df-e^2}} = \nu
$$
we get from \eqref{eq:y_ellipse} that
\begin{equation*}
\COV{F_1,F_2} = 
\frac{16\pi}{\sqrt{4ac-b^2}\sqrt{4df-e^2}}
\sum_{\nu} \frac{r_1(\nu) r_2(\nu)}{\nu^{3}}
\end{equation*}
where $\nu$ goes over all positive common frequencies.
\end{proof}

\section{Covariance in Short Intervals}
In this section we compute the covariance in a short interval $h(T)$ of general domains $\Omega_1,\Omega_2$ as defined in \autoref{defn:covariance_in_short_intervals}.
We will always assume that the global covariance of the normalized error functions $F_1,F_2$ of $\Omega_1,\Omega_2$ is non-zero.
We begin in \autoref{sec:approx} by approximating our functions $F_{1}(t,h),F_{2}(t,h)$ with trigonometric polynomials and calculate their covariance.
Then, in \autoref{sec:gen_cov_short} we use the fact that we found the covariance of the approximating trigonometric polynomials to deduce the covariance of $F_{1}(t,h),F_{2}(t,h)$.

\subsection{Covariance of Approximating Trigonometric Polynomials}\label{sec:approx}

\subsubsection{The Trigonometric Polynomials $P_{i,N}$}

The functions $F_i(t)$, $(i=1,2)$ are approximated by 
$$
P_{i,N}(t) = 
\frac{1}{\pi}\sum_{\substack{n\in\ZZ^2 \\ 0<|n|<N}}\frac{\sqrt{\rho_i(n)}}{|n|^{3/2}}
\cos\left(2\pi Y_i(n) t - \frac{3\pi}{4}\right)
$$
in the sense of \autoref{thm:trig_aprox}.

It will be useful to introduce the following notation:
\begin{definition}[Averaging Operator]
For a function $F$, we define
$$
\BRAK{F}_T = \frac{1}{T}\int_{0}^{T}F(t)dt.
$$
\end{definition}

\autoref{thm:trig_aprox} then states that 
$$
\BRAK{\left(F_i(t)-P_{i,N}(t)\right)^2}_T = \BigO{N^{-1/3}}
$$
provided that $N\leq T^{2/5}$.

Consider the trigonometric polynomials 
$$P_N(t,h) = P_N(t+h) - P_N(t).$$
Using the trigonometric identity 
$$
\cos(\alpha) - \cos(\beta) = 
2\sin\left(\frac{\alpha + \beta}{2}\right)\sin\left(\frac{\beta - \alpha }{2}\right)
$$
we get that
\begin{multline*}
P_{i,N}(t,h) = \\
\frac{2}{\pi}\sum_{\substack{n\in\ZZ^2 \\ 0<|n|<N}}\frac{\sqrt{\rho_i(n)}}{|n|^{3/2}}
\sin\left(2 \pi Y_i(n) \frac{h}{2}\right)
\sin\left(2\pi  Y_i(n) \left( t + \frac{h}{2}\right) + \frac{\pi}{4}\right).
\end{multline*}

From the triangle inequality we get that, for sufficiently small $h$, $P_{i,N}(t,h)$ is an approximation of $F_i(t,h)$ in much the same way as $P_{i,N}(t)$ is an approximation of $F_i(t)$, since
\begin{multline*}
\sqrt{\BRAK{\left(F_i(t,h)-P_{i,N}(t,h)\right)^2}_T}
\leq \\
\sqrt{\BRAK{\left(F_i(t+h)-P_{i,N}(t+h)\right)^2}_T}+
\sqrt{\BRAK{\left(F_i(t)-P_{i,N}(t)\right)^2}_T}
\asymp\\
\sqrt{\BRAK{\left(F_i(t)-P_{i,N}(t)\right)^2}_T}.
\end{multline*}

We continue by computing the covariance of $P_{1,N}(t,h)$ and $P_{2,N}(t,h)$.

\subsubsection{Covariance of $P_{1,N}$ and $P_{2,N}$} \label{sec:aprox_covar}

In this section we calculate 
$$\BRAK{P_{1,N}(t,h) P_{2,N}(t,h)}_T$$
where $N=N(T)$ is sufficiently small.

We introduce the following notation:
$$
D(M) = \min_{\substack{n,m\in\mathbb{Z}^2\\ 0\leq |n|,|m| \leq M \\ Y_1(n)\neq Y_2(m)}}\left|Y_1(n) - Y_2(m)\right|.
$$
\begin{definition}
Let $\Omega_1,\Omega_2$ be two convex planar domains with smooth boundary.
We say that $\Omega_1,\Omega_2$ are a \textit{$\kappa$-Diophantine pair} for some $\kappa > 0$ if $D(M) \gg M^{-\kappa}$.
\end{definition}

\begin{lemma}\label{lem:3sec3lem1}
Assume that $\Omega_1,\Omega_2$ are a $\kappa$-Diophantine pair.
Let $N(T) =\BigO{ T^{\left(\frac{25}{6} + \kappa\right)^{-1}}} $.
Under these assumptions
\begin{multline}
\BRAK{P_{1,N}(t,h(T)) P_{2,N}(t,h(T))}_T = \\
\frac{2}{\pi^2}
\sideset{}{'}\sum_{|n|<N}
\sideset{}{'}\sum_{\substack{|m|<N \\ Y_1(n)=Y_2(m)}}
\frac{\sqrt{\rho_1(n)}}{|n|^{3/2}}
\frac{\sqrt{\rho_2(m)}}{|m|^{3/2}}
\sin\left(2 \pi Y_1(n) \frac{h}{2}\right)^2
+\BigO{N^{-1/6}}
\end{multline}
\end{lemma}

\begin{proof}
Let $T>0$.
From the linearity of $\BRAK{F}_T$ as a function of $F$:
\begin{equation}\label{eq:3lem3eq1}
\begin{multlined}
\BRAK{P_{1,N}(t,h(T)) P_{2,N}(t,h(T))}_T =  \\
\frac{4}{\pi^2} 
\sideset{}{'}\sum_{\substack{n\in\ZZ^2 \\ |n|<N}}\quad
\sideset{}{'}\sum_{\substack{m\in\ZZ^2 \\ |m|<N}}
\frac{\sqrt{\rho_1(n)}}{|n|^{3/2}}
\frac{\sqrt{\rho_2(m)}}{|m|^{3/2}}
\sin\left(2 \pi Y_1(n) \frac{h}{2}\right)
\sin\left(2 \pi Y_2(m) \frac{h}{2}\right) \times \\
\BRAK{\sin\left(2\pi  Y_1(n) \left( t + \frac{h}{2}\right)+\frac{\pi}{4}\right) 
\sin\left(2\pi  Y_2(m) \left( t + \frac{h}{2}\right)+\frac{\pi}{4}\right)}_T
\end{multlined}
\end{equation}

We can rewrite $\sin(\alpha)$ as $\frac{1}{2i}\left(\exp(i\alpha) - \exp(-i\alpha)\right)$.
We get that 
\begin{equation}\label{eq:3lem3eq2}
\begin{split}
\BRAK{\sin\left(2\pi  Y_1(n) 
\left( t + \frac{h}{2}\right)
\right.\right. & \left.\left.
+\frac{\pi}{4}\right)
\sin\left(2\pi  Y_2(m) \left( t + \frac{h}{2}\right)+\frac{\pi}{4}\right)}_T = \\
& \frac{1}{4}\left[
\BRAK{\exp(2\pi i(Y_1(n)-Y_2(m)))}_T e^{i\pi h (Y_1(n)-Y_2(m))} \right.\\
&\quad + \BRAK{\exp(2\pi i(Y_2(m)-Y_1(n)))}_T e^{i\pi h (Y_2(m)-Y_1(n))}\\
&\quad  - \BRAK{\exp(2\pi i(Y_1(n)+Y_2(m)))}_T e^{-i\pi (1/2 + h) (Y_1(n)+Y_2(m))}\\
&\quad  \left. - \BRAK{\exp(-2\pi i(Y_1(n)+Y_2(m)))}_T e^{i\pi (1/2 + h) (Y_1(n)+Y_2(m))}\right] 
\end{split}
\end{equation}
If $Y_1(n) = Y_2(m)$ we get that 
$$
\BRAK{\exp( \pm2\pi i (Y_1(n) - Y_2(m))}_T = 1.
$$
Since both $n,m$ satisfy $|n|,|m|\leq N$,
in all other cases we get-
$$
\left|\pm Y_1(n) \pm Y_2(m)\right| \geq  D(N).
$$
Since $\Omega_1,\Omega_2$ are a $\kappa$-Diophantine pair, and from our choice of $N$, we get that
$$
D(N) \gg N^{-\kappa} \gg T^{-\kappa/\left(\frac{25}{6} + \kappa\right)}.
$$
Thus, since $\left|\BRAK{\exp(i\alpha t)}_T\right| \leq \frac{2}{T\alpha}$, we get that in these cases
$$
\left|\BRAK{\exp(2\pi i(\pm Y_1(n) \pm Y_2(m)))}_T\right| \ll T^{\kappa/\left(\frac{25}{6} + \kappa\right) - 1}
$$

And so, since there are $N^4$ pairs of $n,m$ we get that-
\begin{multline*}
\sideset{}{'}\sum_{\substack{ \pm Y_1(n) \neq \pm Y_2 (m) \\ |n|,|m|\leq N }}
\BRAK{\exp{\left(2\pi i (\pm Y_1(n) \pm Y_2(m))\right)}}_T
\ll
N^4 T^{\kappa/\left(\frac{25}{6} + \kappa\right) - 1}\\
\ll
T^{(4+\kappa)/\left(\frac{25}{6} + \kappa\right) - 1}
= T^{-1/6(\frac{25}{6} + \kappa)}
\ll N^{-1/6}.
\end{multline*}

Thus, from \eqref{eq:3lem3eq1}, \eqref{eq:3lem3eq2} we get that
\begin{multline*}
\BRAK{P_{1,N}(t,h(T)) P_{2,N}(t,h(T))}_T =  \\
\frac{2}{\pi^2}
\sideset{}{'}\sum_{|n|<N}
\sideset{}{'}\sum_{\substack{|m|<N \\ Y_1(n)=Y_2(m)}}
\frac{\sqrt{\rho_1(n)}}{|n|^{3/2}}
\frac{\sqrt{\rho_2(m)}}{|m|^{3/2}}
\sin\left(2 \pi Y_1(n) \frac{h}{2}\right)^2
+\BigO{N^{-1/6}}
\end{multline*}

\end{proof}

Our main result in this section is the following proposition, giving an expression for 
$$
\BRAK{P_{1,N}(t,h(T)) P_{2,N}(t,h(T))}_T
$$
which does not depend on $N$.
There is an implicit dependence, as the formula requires $T$ to be sufficiently larger than $N$, and $N$ to be sufficiently larger than $h$.
\begin{proposition} \label{lem:covar_pol}
Define the function 
$$
f(h) = 
4 \sum_{\nu }\hat{F_1}(\nu)\overline{\hat{F_2}(\nu)}\sin^2\left(\frac{h}{2} \nu\right)
$$
where $\nu$ goes over common frequencies of $F_1$ and $F_2$.

Under the assumptions of \autoref{lem:3sec3lem1} and assuming $h^{-12}(T) \ll N(T)$, 
we have that
\begin{equation*}
\BRAK{P_{1,N}(t,h(T)),P_{2,N}(t,h(T))}_T
= f(h(T)) \left(1 + o(1)\right)
\end{equation*}
\end{proposition}



In order to prove \autoref{lem:covar_pol}, we will first need the following two Lemmas regarding the function $f$-
\begin{lemma}\label{lem:rate_of_decay}
Let $F_1(t),F_2(t)$ be the normalized error functions corresponding to convex planar domains with smooth boundary.
Assume also that $F_1,F_2$ have a non-zero global covariance.
Define 
$$
r(M) = \sum_{|\nu|>M}\hat{F_1}(\nu)\overline{\hat{F_2}}(\nu).
$$
where $\nu$ goes over all common frequencies.
Then
$$
r(M) \ll M^{-1/3}.
$$
\end{lemma}
\begin{proof}
The function $r(M)$ is weakly decreasing, since $\hat{F_1}(\nu)\overline{\hat{F_2}}(\nu)$ is always positive.
Furthermore, from \nameref{thm:parseval} we know that $r(M)$ tends to 0 as $M$ tends to $\infty$.
Using \nameref{thm:parseval}, namely \autoref{coro:parseval}, we get that
$$
r(M) = \sum_{\nu > M}\hat{F_1}(\nu)\overline{\hat{F_2}(\nu)} = 
\lim_{T\rightarrow\infty}\BRAK{\left(F_1(t) - P_{1,M}(t)\right)\left(F_2(t) - P_{2,M}(t)\right)}_T
$$
where we used the fact $P_{1,M}$ and $P_{2,M}$ are the truncated frequency expansions of $F_1$ and $F_2$ respectively, as seen in \eqref{eq:frequency_expansion}.
Using the Cauchy–Schwarz inequality we get that
$$
r(M) \leq \lim_{T\rightarrow\infty}\sqrt{\BRAK{\left(F_1(t) - P_{1,M}(t)\right)^2}_T}\sqrt{\BRAK{\left(F_2(t) - P_{2,M}(t)\right)^2}_T}.
$$
From \autoref{thm:trig_aprox}, we have that
$$
\lim_{T\rightarrow\infty}\BRAK{\left(F_i(t) - P_{i,M}(t)\right)^2}_T = \BigO{M^{-1/3}}.
$$
Thus we have that $r(M)\ll M^{-1/3}$.
\end{proof}

\begin{lemma}\label{lem:lower_bound_f}
Let $F_1(t),F_2(t)$ be the normalized error functions corresponding to convex planar domains with smooth boundary.
Assume also that $F_1,F_2$ have a non-zero global covariance.
Then
$$
\sum_{|\nu|\leq h^{-1}}\hat{F_1}(\nu)\overline{\hat{F_2}(\nu)}\sin^2\left(\frac{h}{2}\nu\right) \gg h^2 \log\left(h^{-1}\right).
$$
\end{lemma}
\begin{proof}
We assume that the global covariance of $F_1$ and $F_2$ is positive.
This means that there exist a common frequency $\nu_0$ such that 
$\hat{F_1}(\nu_0)\overline{\hat{F_2}(\nu_0)} > 0$.
Let $k\in\NN$, and consider the frequency $k\nu_0$.
From the equalities
$$
Y_i(k n) = k Y(n), \quad \rho_i(k n) = \rho_i(n),\quad |kn|^{3/2} = k^{3/2}|n|^{3/2}
$$
we get that $\hat{F_1}(k\nu_0)\overline{\hat{F_2}(k\nu_0)} \geq \frac{1}{k^3}\hat{F_1}(\nu_0)\overline{\hat{F_2}(\nu_0)}$.
This shows that-
\begin{multline*}
f(h) = \sum_{\nu}\hat{F_1}(\nu)\overline{\hat{F_2}(\nu)}\sin^2\left(\frac{h}{2}\nu\right)
\gg
\sum_{ k \leq \nu_0^{-1}h^{-1}}\hat{F_1}(k\nu_0)\overline{\hat{F_2}(k\nu_0)} (h k \nu_0)^2 \\
\gg
h^2  \left(\nu_0^2 \hat{F_1}(\nu_0)\overline{\hat{F_2}(\nu_0)}\right)
\left(\sum_{k \ll h^{-1}}\frac{1}{k}\right) \gg  h^2\log\left(h^{-1}\right)
\end{multline*}

\end{proof}

We can now prove \autoref{lem:covar_pol}:
\begin{proof}[Proof of \autoref{lem:covar_pol}]

From \autoref{lem:3sec3lem1} we have that
\begin{multline*}
\BRAK{P_{1,N}(t,h(T)),P_{2,N}(t,h(T))}_T = \\
\frac{2}{\pi^2}
\sideset{}{'}\sum_{|n|<N}
\sideset{}{'}\sum_{\substack{|m|<N \\ Y_1(n)=Y_2(m)}}
\frac{\sqrt{\rho_1(n)}}{|n|^{3/2}}
\frac{\sqrt{\rho_2(m)}}{|m|^{3/2}}
\sin\left(2 \pi Y_1(n) \frac{h}{2}\right)^2
+\BigO{N^{-1/6}}.
\end{multline*}
We are going to show that $f(h) \gg h^2\log\left(h^{-1}\right)$, and since $N>h^{-12}$ the error term $\BigO{N^{-1/6}}$ is going to be negligible.

We use the fact that $Y_i(n)\asymp |n|$ and rewrite the RHS as
\begin{equation}\label{eq:rhs_freq_exp}
\frac{2}{\pi^2}
\sum_{0< v \leq N}
\sin^2\left(\pi v h\right)
\left(
    \sum_{n: Y_1(n) = v}\frac{\sqrt{\rho_1(n)}}{|n|^{3/2}}
\right)
\left(
    \sum_{m: Y_2(m) = v}\frac{\sqrt{\rho_2(m)}}{|m|^{3/2}}
\right)
\end{equation}
where $v$ goes over all common values of $Y_1(n)$ and $Y_2(m)$.
From \eqref{eq:frequency_expansion} we get that 
$\hat{F_i}(\pm 2\pi  v) = \frac{1}{2\pi}  \sum_{n: Y_i(n) = v}\frac{\sqrt{\rho_i(n)}}{|n|^{3/2}}$.
Thus we can rewrite \eqref{eq:rhs_freq_exp} as:
\begin{multline*}
\frac{2}{\pi^2}
\sum_{0< v \leq N}
\sin^2\left(\pi v h\right)
\frac{4\pi^2}{2}
\left(
    \hat{F_1}(2\pi  v)
    \overline{\hat{F_2}(2\pi  v)}
     + 
     \hat{F_1}(- 2\pi  v)
    \overline{\hat{F_2}(- 2\pi  v)}
\right) = \\
4\sum_{|\nu| \leq 2\pi N}\hat{F_1}(\nu)\overline{\hat{F_2}(\nu)}\sin^2\left(\frac{h}{2} \nu\right)
\end{multline*}
where $\nu$ goes over all common frequencies.

We now wish to prove that this expression is asymptotically equal to $f(h)$.
We do this by showing that the terms in $f$ corresponding to frequencies larger then $N$ are negligible.
And so, it is enough to show that 
$$
\sum_{|\nu| > N}\hat{F_1}(\nu)\overline{\hat{F_2}(\nu)}\sin^2\left(\frac{h}{2} \nu\right)
\lll
\sum_{|\nu| \leq h^{-1}}\hat{F_1}(\nu)\overline{\hat{F_2}(\nu)}\sin^2\left(\frac{h}{2} \nu\right)
$$
where $A\lll B$ means that $\frac{A}{B}\xrightarrow[]{h\rightarrow 0}0$. 
Showing this will prove that 
\begin{multline*}
\BRAK{P_{1,N}(t,h(T)),P_{2,N}(t,h(T))}_T = 4\sum_{|\nu| \leq 2\pi N}\hat{F_1}(\nu)\overline{\hat{F_2}(\nu)}\sin^2\left(\frac{h}{2} \nu\right) \\
= f(h) - \sum_{|\nu| > 2\pi N}\hat{F_1}(\nu)\overline{\hat{F_2}(\nu)}\sin^2\left(\frac{h}{2} \nu\right) \sim f(h).
\end{multline*}

From \autoref{lem:lower_bound_f} we have that 
$$
\sum_{|\nu| \leq h^{-1}}\hat{F_1}(\nu)\overline{\hat{F_2}(\nu)}\sin^2\left(\frac{h}{2} \nu\right) \gg h^2\log\left(h^{-1}\right).
$$

Using the bound $\left|\sin(x)\right|\leq 1$ and from \autoref{lem:rate_of_decay} we get
$$
\sum_{|\nu| > N}\hat{F_1}(\nu)\overline{\hat{F_2}(\nu)}\sin^2\left(\frac{h}{2} \nu\right) \ll r(N) \ll N^{-1/3}.
$$

Since $N>h^{-12}$ we have $N^{-1/3} \lll h^2\log\left(h^{-1}\right)$ which completes the proof.
\end{proof}

\begin{corollary}\label{cor:minimal_f}
The function $f(h)$ defined in \autoref{lem:covar_pol} satisfies
$$f(h) \gg h^2\log\left(h^{-1}\right).$$
Furthermore, if we assume that all common frequencies of $F_1,F_2$ are of the form $k\nu_0$ for some initial $\nu_0$ and $k\in\ZZ$,
and we assume that $2\pi Y_1(n) = 2\pi Y_2(m) = k\nu_0$ only when $n=k n_0, m = k m_0$ for some solution $n_0, m_0$ of $2\pi Y_1(n_0)=2\pi Y_2(m_0) = \nu_0$ (i.e. the number of solutions is constant), then:
$$
f(h) = C h^2 \log\left(h^{-1}\right) + \BigO{h^2}
$$
where 
$$
C = 2\nu_0^2 \hat{F_1}(\nu_0)\overline{\hat{F_2}(\nu_0)}.
$$

\end{corollary}
\begin{proof}
The bound $f(h) \gg h^2\log\left(h^{-1}\right)$ follows from \autoref{lem:lower_bound_f}.

As for the second part of the statement, 
the fact that 
$$2\pi Y_1(n) =2\pi Y_2(m) = k\nu_0$$
only when $n=k n_0, m = k m_0$ for some solution $n_0, m_0$ of $2\pi Y_1(n_0)=2\pi Y_2(m_0) = \nu_0$ implies that
$$
\hat{F_1}(k\nu_0) = \frac{1}{k^3}\hat{F_1}(\nu_0),\quad
\hat{F_2}(k\nu_0) = \frac{1}{k^3}\hat{F_2}(\nu_0).
$$
We then have that-
\begin{multline*}
f(h) = 
4\sideset{}{'}\sum_{k\in\ZZ}
\hat{F_1}(k\nu_0)\overline{\hat{F_2}(k\nu_0)}
\sin^2\left(\frac{h}{2}k\nu_0\right) 
= \\
8 \hat{F_1}(\nu_0)\overline{\hat{F_2}(\nu_0)} 
\sum_{k\in\NN}
\frac{1}{k^3}\sin^2\left(\frac{h}{2}k\nu_0\right) 
=\\
C h^2\sum_{k \ll h^{-1}}\frac{1}{k} + 
\BigO{\sum_{k \gg h^{-1}} \frac{1}{k^3}} + 
\BigO{h^2} 
=
C h^2 \log\left( h^{-1}\right) + \BigO{h^2} .
\end{multline*}
\end{proof}

\subsection{Covariance of the Error Term in Short Intervals of General Domains}\label{sec:gen_cov_short}
We now compute the covariance of $F_1(t,h(T)),F_2(t,h(T))$.

\begin{theorem}\label{thm:covar_short_interval_F1_F2}
Let $\Omega_1,\Omega_2$ be two convex planar domains, each with a boundary which is smooth and has positive curvature and let 
$F_1(t)$, $F_2(t)$ be their normalized error functions.
Assume that the global covariance of $F_1(t)$, $F_2(t)$ is non-zero.
Assume also that $\Omega_1,\Omega_2$ are a $\kappa$-Diophantine pair.
Let $h(T)$ be a function of $T$ tending to 0 such that $h(T) \gg T^{-1/(50 + 12\kappa)}$.

Then-
$$
\COV{F_1(T,h(T)),F_2(T,h(T))} = f(h)
$$
where 
$$
f(h) = 4\sum_{\nu}\hat{F_1}(\nu)\overline{\hat{F_2}(\nu)}\sin^2\left( \frac{h}{2} \nu\right).
$$

\end{theorem}
\begin{proof}
Take $N(T) = h^{-12}(T) \ll T^{1/(\frac{25}{6}+\kappa)}$.
Note that the conditions on $h$ imply that $N(T) \ll T^{2/5}$.
From \autoref{lem:covar_pol} we have
$$
\BRAK{P_{1,N}(t,h(T))P_{2,N}(t,h(T))}_T = f(h)(1+o(1)).
$$

Thus, in order to prove the theorem it is enough to show that-
$$
\lim_{T\rightarrow\infty}\frac{1}{f(h(T))}\left(\BRAK{F_1(t,h)F_2(t,h)}_T - \BRAK{P_{1,N}(t,h)P_2(t,h)}_T\right) = 0.
$$

From the equality
\begin{equation*}
P_{1,N}P_{2,N} - F_1 F_2   = 
(F_1-P_{1,N})(F_2-P_{2,N}) - (F_1-P_{1,N})F_2 - (F_2-P_{2,N})F_1.
\end{equation*}
We get-
\begin{multline*}
\frac{1}{f(h(T))}\left|\BRAK{F_1(t,h)F_2(t,h)}_T - \BRAK{P_{1,N}(t,h)F_2(t,h)}_T\right| 
\leq \\
\frac{1}{f(h(T))}\left|
\BRAK{(F_1-P_{1,N})(F_2-P_{2,N})}_T
+\right.\\
\left. \BRAK{(F_1-P_{1,N})F_2}_T
+
\BRAK{(F_2-P_{2,N})F_2}_T \right|.
\end{multline*}

Using the Cauchy–Schwarz inequality we get that
\begin{align}
\begin{split} \label{eq:triangle_cs_ineq}
\frac{1}{f(h(T))}\left|
\right.&\left.
\BRAK{F_1(t,h)F_2(t,h)}_T- 
\BRAK{P_{1,N}(t,h)F_2(t,h)}_T\right| 
\leq \\
& \frac{1}{f(h)}\left(
\sqrt{\BRAK{(F_1-P_{1,N})^2}_T}
\sqrt{\BRAK{(F_2-P_{2,N})^2}_T} \right.\\
& \qquad\quad +
\sqrt{\BRAK{(F_1-P_{1,N})^2}_T}
\sqrt{\BRAK{F_2^2(t,h)}_T} \\
& \left. \qquad\quad + 
\sqrt{\BRAK{(F_2-P_{2,N})^2}_T}
\sqrt{\BRAK{F_1^2(t,h)}_T} \quad
\right).
\end{split}
\end{align}

From \autoref{thm:trig_aprox} we get that $\BRAK{(F_i-P_{i,N})^2}_T = \BigO{N^{-1/3}}$.
As for the terms $\sqrt{\BRAK{(F_i^2(t,h)}_T}$, using the triangle inequality we get
\begin{multline*}
\sqrt{\BRAK{(F_i^2(t,h)}_T} = 
\sqrt{\BRAK{(F_i(t+h) - F_i(t))^2}_T} \leq \\
\sqrt{\BRAK{F_i^2(t+h)}_T} + \sqrt{\BRAK{F_i^2(t)}_T} \asymp
\sqrt{\BRAK{F_i^2(t)}_T}
\end{multline*}
which is bounded in $T$ since $F_i$ have finite variance (see \cite{bleher}).

And so, from \eqref{eq:triangle_cs_ineq} we get that
\begin{multline*}
\frac{1}{f(h(T))}\left(\BRAK{F_1(t,h)F_2(t,h)}_T- 
 \BRAK{P_{1,N}(t,h)F_2(t,h)}_T\right)
 = \\
 \frac{1}{f(h)}\left(\BigO{N^{-1/3}} + \BigO{N^{-1/6}} + \BigO{N^{-1/6}}\right).
\end{multline*}
From \autoref{lem:lower_bound_f} we get that $f(h)\gg h^2\log\left(h^{-1}\right)$.
And so, 
\begin{multline*}
\frac{1}{f(h(T))}\left(\BRAK{F_1(t,h)F_2(t,h)}_T- 
\BRAK{P_{1,N}(t,h)F_2(t,h)}_T\right)
= \\
\BigO{\frac{N^{-1/6}}{h^2\log\left(h^{-1}\right)}} = o(1)
\end{multline*}
where the last equality used the fact that $N = h^{-12}$.
\end{proof}

\section{Covariance of Ellipses in Short Intervals}\label{chap:ellipse_short}

In this section we consider the covariance in short intervals of normalized error functions of ellipses.
In \autoref{sec:gen_elips_short_intervals} we prove a result concerning the covariance of 'generic' ellipses.
In \autoref{sec:rec_triangle_covar} we consider a more specific case.
Namely, we consider the covariance between $x^2+\alpha y^2 \leq 1$ and $x^2 + xy + y^2 \leq \frac{3}{4}$ for various positive $\alpha$'s.
This example is of particular interest since it related to error terms of the Dirichlet eigenvalue problem for a rectangle and for an equilateral triangle.

\subsection{Covariance of Generic Ellipses}\label{sec:gen_elips_short_intervals}
Our main goal in this section is to prove \autoref{thm:ellipse_short_covar} regarding the covariance of a generic ellipse pair.
By this, we mean that the set of ellipses that don't satisfy the desired properties has measure 0.
In order to make this statement precise, we will first need to define the measure space that we are referring too.
We then proceed by proving that for $\kappa>5$ almost all ellipse pairs are $\kappa$-Diophantine, which allows us to invoke \autoref{thm:covar_short_interval_F1_F2} and get a formula for the covariance in short intervals (with some restriction on the decay rate of the short interval).
We also prove that almost all pairs of ellipses satisfy the conditions of \autoref{cor:minimal_f}, proving that for almost all pairs the covariance has order $h^2\log\left(h^{-1}\right)$.
\subsubsection{Definition of the Measure Space}
Let 
$$
\Omega_1: a x^2+ b xy + c\leq 1, \quad
\Omega_2: d x^2 + e xy + f y^2\leq 1
$$

be two ellipses for some parameters $a,b,c,d,e,f$.

\begin{definition}
We define $\Gamma$ as the set of ellipse pairs.
We identify $\Gamma$ with a subset of $\RR^6$ (with its standard Lebesgue measure) defined by
$$
\Gamma = \left\{
(a,b,c,d,e,f) \in \RR^6 \mid
4ac > b^2, 4de > f^2
\right\}
$$
where we identify $(a,b,c,d,e,f)$ with $(\Omega_1,\Omega_2)$.
\end{definition}

We will however be interested in ellipse pairs that have non-zero global covariance.
To this end, for any non-zero $n=(n_1,n_2), m=(m_1,m_2)\in\ZZ^2$ we define 
$\Gamma_{(n,m)}\subset \Gamma$ as the zero set of 
\begin{multline*}
\left(Y_1(n) - Y_2(m)\right)(a,b,c,d,e,f)
= \\
\frac{2\sqrt{a n_2^2 - b n_1 n_2 + c n_1^2}}{\sqrt{4ac-b^2}}
-
\frac{2\sqrt{d m_2^2 - e m_1 m_2 + f m_1^2}}{\sqrt{4df-e^2}}
\end{multline*}
where we used \eqref{eq:y_ellipse} for the formulas of $Y_1$ and $Y_2$.

\begin{definition}
We define the metric space of ellipse pairs with non-zero global covariance as 
\begin{equation}\label{eq:gamma_prime}
\Gamma' = \bigcup_{\substack{n,m\in\ZZ^2 \\ n,m\neq 0}}\Gamma_{(n,m)}
\end{equation}
with the inherited metric from $\Gamma$.
\end{definition}

\subsubsection{Diophantine Property for Ellipse Pairs}

\begin{proposition}\label{prop:gen_ellipsse_diophantine}
For all $\kappa > 5$, almost all ellipse pairs that have non-zero global covariance are a $\kappa$-Diophantine pair.
\end{proposition}

We will first need the following lemma.

Define a function $\psi:\ZZ^4\rightarrow\ZZ^6$ by
\begin{equation}\label{eq:psi_def}
\psi(n_1,n_2,m_1,m_2) = 
(n_2^2,n_1n_2,n_1^2, m_2^2, m_1m_2, m_2^2).
\end{equation}

\begin{definition}
We say that $\mathbf{p}\in\RR^6$ is V.W.A (very well approximable) if there is some $\epsilon>0$ such that there are infinitely many $(n_1,n_2,m_1,m_2)\in\ZZ^4$ for which 
$$
0< \left|
\mathbf{p}\cdot \psi(n_1,n_2,m_1,m_2)
\right|
\leq \left(\max\left(n_1,n_2,m_1,m_2\right)\right)^{-4 - \epsilon}.
$$
\end{definition}

\begin{lemma}\label{lem:hyperplane_vwa}
Let $\mathbf{q}\neq 0\in\ZZ^6$ and let $H = \mathbf{q}^{\perp}$ be the hyper-plane  perpendicular to $\mathbf{q}$.
Almost all $\mathbf{p}\in H$ are not V.W.A.
\end{lemma}

\begin{proof}
It is enough to prove the statement for some local compact subset of $H$, which we denote $\mathcal{C}$.
Let $n=(n_1,n_2),m=(m_1,m_2)\in\ZZ^2$.
Denote by $\mathcal{C}_{n,m}$ the set of points $\mathbf{p}\in \mathcal{C}$ which satisfy
$$
0< \left|
\mathbf{p}\cdot \psi(n,m)
\right|
\leq \left(\max\left(n_1,n_2,m_1,m_2\right)\right)^{-4 - \epsilon}.
$$
If $\psi(n,m)$ is proportional to $\mathbf{q}$, then $\mathcal{C}_{n,m}=\emptyset$ has measure 0.
Otherwise, the points in $\mathcal{C}_{n,m}$ satisfy
$$
0<\left| \mathbf{p}\cdot \frac{\psi(n,m)}{\|\psi(n,m)\|}
\right|
\ll \delta
$$
where $\delta = \left(\max\left(n_1,n_2,m_1,m_2\right)\right)^{-6 - \epsilon}$.
This means that $\mathbf{p}$ is in a $\delta$ neighbourhood of the hyper-plane $\psi(n,m)^{\perp}$ in $\RR^6$. 
We now show that this also means that $\mathbf{p}$ is in a small neighbourhood of $H\cap\psi(n,m)^{\perp}$ in $H$.

Denote by $\theta$ the angle between $\mathbf{q}$ and $\psi(n,m)$.
Since $\mathbf{q}\cdot\psi(n,m)$ is an integer, and $\mathbf{q}$ is not proportional to $\psi(n,m)$ we get that
\begin{equation}\label{eq:theta_bound}
\cos^2\left(\theta\right) =
\frac{\left(\mathbf{q}\cdot\psi(n,m)\right)^2}{\|\mathbf{q}\|^2 \| \psi(n,m)\|^2} 
\leq 
\frac{\|\mathbf{q}\|^2 \|\psi(n,m)\|^2 -1}{\|\mathbf{q}\|^2 \| \psi(n,m)\|^2}
\leq
1 - \frac{1}{\|\mathbf{q}\|^2 \|\psi(n,m)\|^2}.
\end{equation}
It can be seen that any point on $H$ that is in a $\delta$ neighbourhood of $\psi(n,m)^{\perp}$, is in a $\delta/\sin\left(\theta\right)$ neighbourhood of $H\cap \psi(n,m)^{\perp}$ (see \autoref{fig:hyperplane_disc}). 
\begin{figure}[hbt]
\includegraphics[trim=0 75 0 0, width=0.75\textwidth]{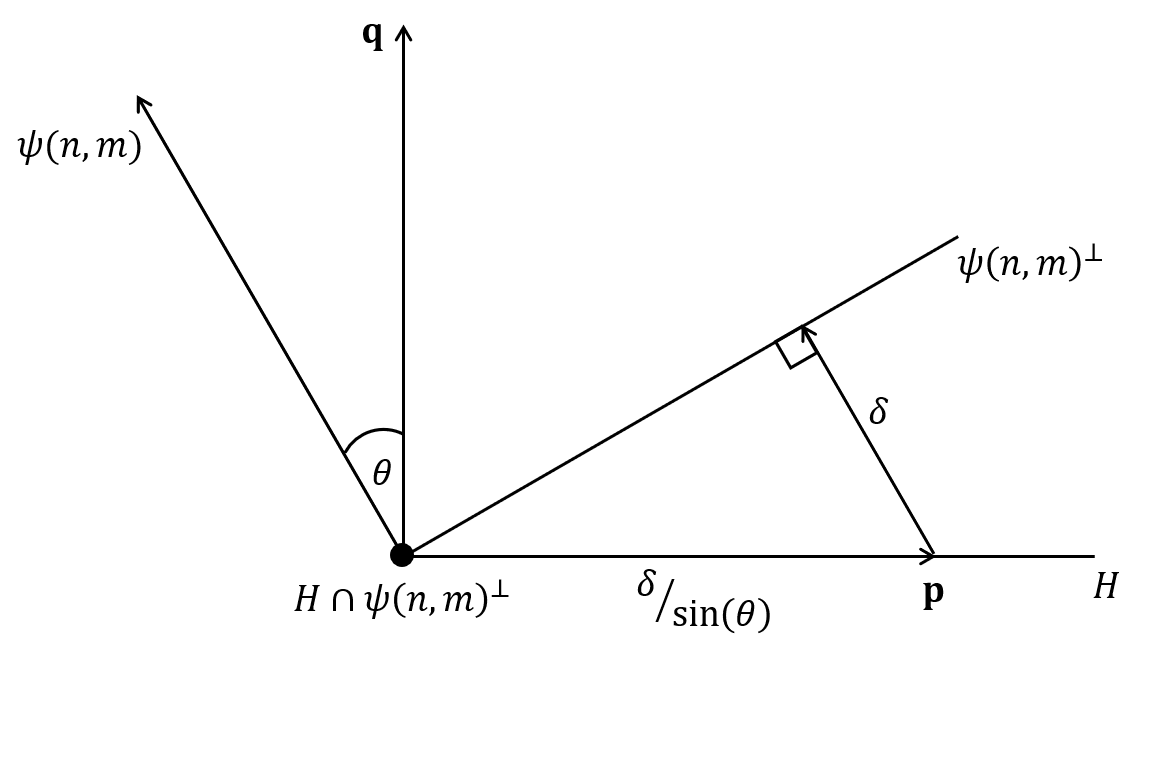}
\centering
\caption{A projection on to $\text{span}\left(\mathbf{q},\psi(n,m)\right)$}
\label{fig:hyperplane_disc}
\end{figure}

Thus, from \eqref{eq:theta_bound} we get that $\mathbf{p}$ is in a 
$$
\delta / \sin\left(\theta\right) \asymp_{\mathbf{q}}
 \delta \|\psi(n,m)\| \asymp  
\left(\max\left(n_1,n_2,m_1,m_2\right)\right)^{-4 - \epsilon}
$$ 
neighbourhood of $H\cap \psi(n,m)^{\perp}$. 

Since we restricted our attention to a compact subset of $H$, it follows that 
$$
\mu\left(\mathcal{C}_{n,m}\right) \ll_{\mathcal{C}} 
\left(\max\left(n_1,n_2,m_1,m_2\right)\right)^{-4 - \epsilon}
$$
where $\mu$ is the standard Lebesgue measure on $H$. 
We now have that
$$
\sum_{\substack{n,m\in\ZZ^2\\ n,m\neq 0}}
\mu\left(\mathcal{C}_{n,m}\right)
\ll
\sum_{\substack{n,m\in\ZZ^2\\ n,m\neq 0}} \left(\max\left(n_1,n_2,m_1,m_2\right)\right)^{-4 - \epsilon} < \infty.
$$
The statement then follows from the Borel–Cantelli lemma.
\end{proof}

\def\INI#1{{\mathbf{#1}^s}}
\begin{proof}[Proof of \autoref{prop:gen_ellipsse_diophantine}]
Denote by $\mathcal{A}$ the set of points in $\Gamma'$ which are not a $\kappa$-Diophantine pair.
We will show that the set $\mathcal{A}$ has measure zero.
Since \eqref{eq:gamma_prime} expresses $\Gamma'$ as a countable union, it is enough to prove the statement for each set in the union.
Thus, we set some non-zero $n^s=(n^s_1,n^s_2), m^s=(m^s_1,m^s_2)\in\ZZ^2$ and we wish to show that that $\mathcal{A}\cap \Gamma_{(n^s,m^s)}$ has measure 0.

In order for a point to lie in $\mathcal{A}$, we must have that for infinitely many $n=(n_1,n_2),m=(m_1,m_2) \in \ZZ^2$
$$
Y_1(n) \neq Y_2(m),\quad 
\left|Y_1(n) - Y_2(m)\right| \ll_{a,b,c,d,e,f} \left(\max(n_1,n_2,m_1,m_2)\right)^{-\kappa}.
$$

Since $Y_1(n) \asymp_{a,b,c,d,e,f} |n|$, we get that this is equivalent to having infinitely many $n,m \in \ZZ^2$ with
\begin{multline}\label{eq:4prop1eq1}
Y_1^2(n) - Y_2^2(m) = \\
\left(Y_1(n) - Y_2(m)\right)
\left(Y_1(n) + Y_2(m)\right) \ll_{a,b,c,d,e,f} \left(\max\{(n_1,n_2,m_1,m_2\}\right)^{-\kappa + 1}.
\end{multline}

We define the functions $A,B,C,D,E,F$ of $a,b,c,d,e,f$ by
\begin{align*}
A = \frac{a}{4ac-b^2},\quad 
B = \frac{-b}{4ac-b^2},\quad
C = \frac{c}{4ac-b^2}, \\ 
D = \frac{-d}{4df-e^2}, \quad
E = \frac{e}{4df-e^2}, \quad
F = \frac{-f}{4df-e^2}. 
\end{align*}
With these notations we have that
$$
Y_1^2(n) - Y_2^2(m) = (A,B,C,D,E,F)\cdot \psi(n,m).
$$

Consider the function $\Phi:\Gamma \rightarrow \RR^6$ defined by $\Phi:(a,b,c,d,e,f)\mapsto(A,B,C,D,E,F)$.
It is nowhere degenerate since
$$
\left|\frac{\partial(A,B,C,D,E,F)}{\partial(a,b,c,d,e,f)}\right| = \frac{-1}{\left(4ac-b^2\right)^{3}\left(4df-e^2\right)^{3}}<0.
$$

The restriction of $\Phi$ to $\Gamma_{(n^s,m^s)}$ takes elements from the zero set of 
$$Y_1^2(n^s) - Y_2^2(m^s)$$
to the hyper-plane
$$
H : 
\left\{(A,B,C,D,E,F) \mid
(A,B,C,D,E,F) \cdot \psi(n^s,m^s) = 0\right\}.
$$

From \autoref{lem:hyperplane_vwa} we get that almost all points $(A,B,C,D,E,F)\in H$ are not V.W.A.
Thus, since $\Phi$ is nowhere degenerate, it follows that for almost all $(a,b,c,d,e,f)\in\Gamma_{(n^s,m^s)}$, the corresponding $(A,B,C,D,E,F)\in H$ is not V.W.A.
Since $\kappa > 5$, this means that the set of $(a,b,c,d,e,f)\in\Gamma_{(n^s,m^s)}$ that satisfy \eqref{eq:4prop1eq1} for infinitely many $(n,m)$ has measure 0.
It follows that almost all ellipse pairs $(a,b,c,d,e,f)\in\Gamma_{(n^s,m^s)}$ are a $\kappa$-Diophantine pair.

\end{proof}

\subsubsection{Common Frequencies of Generic Ellipses}
In this section we prove a result regarding common frequencies of generic ellipses that have non-zero global covariance.
We will use the same notations we used in the proof of \autoref{prop:gen_ellipsse_diophantine}.

\begin{proposition}\label{prop:common_freq_gen_ellipse}
Almost all ellipse pairs in $\Gamma'$ satisfy the following:
The only $n,m\in\ZZ^2$ for which $Y_1(n) = Y_2(m)$ are $n=k n_0,m=\pm k m_0$ for $k\in\ZZ$ and some initial $n_0,m_0\in\ZZ^2$.
\end{proposition}

We will first note the following-
\begin{lemma}\label{lem:parallel_hyper_planes}
Let $\psi$ be the function defined in \eqref{eq:psi_def}.
Then $\psi(n_1,n_2,m_1,m_2)$ is proportional to $\psi(\ell_1,\ell_2,k_1,k_2)$ if and only if $(n_1,n_2,m_1,m_2)$ is proportional to either $(\ell_1,\ell_2,k_1,k_2)$ or to $(\ell_1,\ell_2,-k_1,-k_2)$.
\end{lemma}
\begin{proof}
This follows from the fact that $(x,y)\mapsto (y^2,xy,x^2)$ is invertible up to $\pm$ sign.
\end{proof}

\begin{proof}[Proof of \autoref{prop:common_freq_gen_ellipse}]
As in the proof of \autoref{prop:gen_ellipsse_diophantine} it is enough to prove this for the case of $\Gamma_{(n^s,m^s)}$ for some $n^s,m^s\in\ZZ^2$.
We wish to show that for almost all $(a,b,c,d,e,f)\in \Gamma_{(n^s,m^s)}$, the only $n,m\in\ZZ^2$ for which $Y_1(n) = Y_2(m)$ are integer multiples of $n^s,m^s$ or of $n^s,-m^s$.
In other words, we wish to show that
\begin{equation}\label{eq:set_a_union}
\mathcal{A} = \bigcup_{\substack{n,m\in\ZZ^2 \\ n,m \neq kn_0,\pm km_0}} \left\{ F_{(n,m)}(a,b,c,d,e,f) = 0\right\}
\end{equation}
has measure 0 in $\Gamma_{(n^s,m^s)}$ where
\begin{multline*}
F_{(n,m)}(a,b,c,d,e,f) = 
Y_1^2(n) - Y_2^2(m) =  \\
\frac{a}{4ac-b^2}n_2^2 - 
\frac{b}{4ac-b^2}n_1n_2 +
\frac{c}{4ac-b^2}n_2^2  \\
- \frac{d}{4df-e^2}m_2^2 + 
\frac{e}{4df-e^2}m_1 m_2 -
\frac{f}{4df-e^2}m_1^2.
\end{multline*}

The expression \eqref{eq:set_a_union} expresses $\mathcal{A}$ as a countable union, so once more it is enough to show that 
$$
\left\{ F_{(n,m)}(a,b,c,d,e,f)=0\right\}
\cap\Gamma_{(n^s,m^s)}
$$
has measure 0 for some specific $n,m$ which is not a multiple of $(n^s,\pm m^s)$.
We once again use the map $\Phi$ defined in \autoref{prop:gen_ellipsse_diophantine}.
The map $\Phi$ takes $\left\{ F_{(n,m)}(a,b,c,d,e,f)=0\right\}$ to the hyper-plane $H_1$ perpendicular to $\psi(n,m)$.
However, since we are restricting our attention to $\Gamma_{(n^s,m^s)}$, the image is also contained in the hyper-plane $H_2$ perpendicular to $\psi(n^s,m^s)$.
Since $n,m$ is not a multiple of $n^s,\pm m^s$, we get from \autoref{lem:parallel_hyper_planes} that $H_1$ and $H_2$ are not equal.
This means that $H_1\cap H_2$ has co-dimension 1 in $H_2$.
From the non-degeneracy of $\Phi$, this means that $\left\{ F_{(n,m)}(a,b,c,d,e,f)=0\right\}\cap \Gamma_{(n^s,m^s)}$ is a sub-manifold of co-dimension 1 in $\Gamma_{(n^s,m^s)}$.
Thus, it has measure 0 as required.

\end{proof}

\subsubsection{Covariance of a Generic Ellipse Pair}
We now prove the main result of this section

\begin{proof}[Proof of \autoref{thm:ellipse_short_covar}]
Let $\epsilon > 0$.
From \autoref{prop:gen_ellipsse_diophantine} we have that almost all ellipse pairs $(\Omega_1,\Omega_2)\in\Gamma'$ are a $(5+\epsilon)$-Diophantine pair.
Thus, from \autoref{thm:covar_short_interval_F1_F2} we get if $h(T) \gg T^{-1/110 + \epsilon'}$, then
$$
\COV{F_1(t,h),F_2(t,h)} = 
4\sum_{\nu}\hat{F_1}(\nu)\overline{\hat{F_2}}(\nu)
\sin^2\left(\frac{h}{2}\nu\right).
$$
Denote by $\nu_0$ the smallest positive common frequency of $F_1$ and $F_2$.
From \autoref{prop:common_freq_gen_ellipse} we get that almost all ellipse pairs satisfy the conditions of \autoref{cor:minimal_f}.
Namely, we get that  if $2\pi Y_1(\pm n_0) = 2\pi Y_2(\pm m_0) = \nu_0$, then the only solutions $n,m$ to $2\pi Y_1(n) = 2\pi Y_2(m)= k\nu_0$ are $n=\pm k n_0, m = \pm k m_0$.
Thus, from \autoref{cor:minimal_f} we have that
$$
4\sum_{\nu}\hat{F_1}(\nu)\overline{\hat{F_2}}(\nu)
\sin^2\left(\frac{h}{2}\nu\right)
= 
C h^2 \log\left( h^{-1}\right) + \BigO{h^2}
$$
where $C$ is given by $2\nu_0^2 \hat{F_1}(\nu_0)\overline{\hat{F_2}(\nu_0)}$.
From \eqref{eq:frequency_expansion},\eqref{eq:y_ellipse},\eqref{eq:ellipse_freq_coef} we get that 
$$C = \frac{16\pi}{\nu_0\sqrt{4ac-b^2}\sqrt{4df-e^2}}.$$
\end{proof}

\subsection{Covariance Between a Rectangle and an Equilateral Triangle}\label{sec:rec_triangle_covar}

We now consider the covariance between the normalized error functions $F_{1},F_{2}$ of the following two  ellipses:
$$\Omega_1: x^2 + xy + y^2 \leq \frac{3}{4}, \quad \Omega_2 : x^2 + \alpha y^2 \leq 1.$$
These normalized error functions are related to the Dirichlet eigenvalue error functions corresponding to an equilateral triangle of side length $\frac{2\pi}{\sqrt{3}}$, and a rectangle with side lengths $\pi$ and $\sqrt{\alpha}\pi$.
The covariance can behave differently depending on $\alpha$.
Our main result from this section is \autoref{thm:specific_ellipse}.
We begin by examining the 'generic' case, and then proceed to explore the case of rational $\alpha$.

\subsubsection{The Generic Case}\label{sec:alpha_irrational}

\begin{lemma}\label{lem:sec4sub2_irrational}
For any $\kappa > 3$, for almost all $\alpha\in\RR$, $\Omega_1,\Omega_2$ is a $\kappa$-Diophantine pair.
\end{lemma}

\begin{proof}
Assume that $\alpha\notin \mathbb{Q}$.
Denote $K = (\kappa -1) /2$.
Since $K > 1$, it is known that almost all $\alpha$ satisfy
\begin{equation}\label{eq:diophantine_alpha_real}
\left| k - \ell \alpha \right| \gg \frac{1}{\ell^K}    
\end{equation}
for almost all $k,\ell \in \ZZ$.
We show that if $\alpha$ is such a number, then $\Omega_1,\Omega_2$ is a $\kappa$-Diophantine pair.

From \eqref{eq:y_ellipse} we get that for $n=(n_1,n_2),m=(m_1,m_2)\in\ZZ^2$:
$$
Y_1(n) = \sqrt{ n_1^2 -  n_1 n_2  + n_2^2}, \quad
Y_2(m) = \frac{1}{\sqrt{\alpha}}\sqrt{\alpha m_1^2 + m_2^2}.
$$
If $n_1,n_2,m_1,m_2\leq M$ then we have that-
\begin{multline}\label{eq:sec4sub2eq1}
\left|Y_1(n) - Y_2(m) \right|= 
\frac{1}{Y_1(n) + Y_2(m)}
\left| n_1^2 -  n_1 n_2 +  n_2^2 - m_2^2 - \alpha m_1^2 \right| 
\\ \gg 
\frac{1}{M}\left| n_1^2 -  n_1 n_2 +  n_2^2 - m_2^2 - \alpha m_1^2 \right|.
\end{multline}
We wish to prove that if $Y_1(n)\neq Y_2(m)$ then $\left|Y_1(n)- Y_2(m)\right| \gg  M^{-\kappa}$.
Since $n_1^2 -  n_1 n_2 +  n_2^2 - m_2^2$ and $m_1^2$ are integers of size at most $M^2$, and since we assume that $\alpha$ satisfies \eqref{eq:diophantine_alpha_real}, it follows that 
$$
n_1^2 -  n_1 n_2  + n_2^2 - m_2^2 - \alpha m_1^2 \gg \left({M^2}\right)^{-K} = M^{-2K}.
$$
From \eqref{eq:sec4sub2eq1} we get that
$$
\left|Y_1(n)- Y_2(m)\right| \gg M^{-(2K+1)} = M^{-\kappa}
$$
which implies that $(\Omega_1,\Omega_2)$ is a $\kappa$-Diophantine pair.
\end{proof}

\begin{proof}[Proof of \autoref{thm:specific_ellipse}, Part 1]
From \autoref{lem:sec4sub2_irrational} we know that almost all $\alpha\not\in \mathbb{Q}$  will result in  $(\Omega_1,\Omega_2)$ being a $\kappa$-Diophantine pair for any $\kappa > 3$.
Let $\alpha$ be such a number, and let $h(T) \gg T^{-1/86 + \epsilon}$ for some $\epsilon>0$.
We can find $\kappa$ sufficiently close to 3 such that 
$$
h(T) \gg T^{-1/86 + \epsilon} \gg T^{-(50 + 12 \kappa)^{-1}}.
$$
Thus, from \autoref{thm:covar_short_interval_F1_F2} we get that 
$$
\COV{F_1(t,h),F_2(t,h)} = 4\sum_{\nu}\hat{F_1}(\nu)\overline{\hat{F_2}}(\nu)
\sin^2\left(\frac{h}{2}\nu\right).
$$
Since $\alpha\notin \mathbb{Q}$, we get that the only common frequencies of $F_1$ and $F_2$ are of the form $2\pi k, k\in\ZZ$.

From \eqref{eq:frequency_expansion},\eqref{eq:ellipse_freq_coef} we get that for $k\in\NN$: 

\begin{equation*}
\hat{F_1}(\pm 2\pi k) = 
e^{\mp 3\pi i/4}
\frac{\sqrt{8\pi}}{\frac{4}{\sqrt{3}}}
\frac{r_1(2\pi k)}{(2\pi k)^{3/2}}
=
e^{\mp 3\pi i/4}
\frac{\sqrt{3}}{4\pi}
\frac{r_1(2\pi k)}{k^{3/2}}
\end{equation*}

where 
$$
r_1(\nu) = \#\left\{ n\in \ZZ^2 \mid 2\pi\sqrt{n_1^2 + n_1 n_2 + n_2^2} = \nu\right\}.
$$

Similarly, noting that there are exactly two solutions $m\in\ZZ^2$ to $Y_2(m) = k$, we get that

\begin{equation*}
\hat{F_2}(\pm 2\pi k) = 
e^{\mp 3\pi i/4}
\frac{\sqrt{8\pi}}{2\sqrt{\alpha}}
\frac{2}{(2\pi k)^{3/2}}
=
e^{\mp 3\pi i/4}
\frac{1}{\pi \sqrt{\alpha} k^{3/2}}.
\end{equation*}

Thus, the covariance is given by
\begin{multline}\label{eq:sec4sub2eq1_cov_alpha}
\COV{F_1(t,h),F_2(t,h)}=\\
4\sum_{\nu}
\hat{F_1}(\nu)\overline{\hat{F_2}}(\nu)
\sin^2\left(\frac{h}{2}\nu\right)
= 
\frac{2\sqrt{3}}{\pi^2\sqrt{\alpha}}
\sum_{k\in\NN} \frac{r_1(2\pi k)}{k^3}
\sin^2\left(\pi h k\right).
\end{multline}

We can further estimate this sum by looking at the first $h^{-1}$ terms and the rest of the terms separately.
We will use the fact that
$$
\sum_{k\leq N}r_1(2\pi k) = \frac{3\sqrt{3}}{\pi} N\log N + \BigO{N}
$$
which we prove in \autoref{thm:mult_case}.
For the terms with $k\gg h^{-1}$ we use the inequality $\sin(x)\leq 1$ and summation by parts to get
\begin{equation}\label{eq:estimate_last_L_terms}
\sum_{k\gg h^{-1}} \frac{r_1(2\pi k)}{k^3}
\sin^2\left(\pi h k\right)
=
\BigO{h^2\log\left(h^{-1}\right)}.
\end{equation}

As for the first $h^{-1}$ terms, using the approximation $\sin(x) = x +\BigO{x^2}$, and once more using summation by parts, we get
\begin{multline}\label{eq:estimate_first_L_terms}
\sum_{k\ll h^{-1}} \frac{r_1(2\pi k)}{k^3}
\sin^2\left(\pi h k\right)
=
3\sqrt{3}\pi h^2 
\sum_{k\ll h^{-1}} \frac{\log(k)}{k} +
\BigO{h^2\log\left(h^{-1}\right)}
= \\
\frac{3\sqrt{3}\pi}{2} h^2 \log^2\left(h^{-1}\right)
 + \BigO{h^2\log\left(h^{-1}\right)}.
\end{multline}

From \eqref{eq:sec4sub2eq1_cov_alpha}, \eqref{eq:estimate_last_L_terms}, \eqref{eq:estimate_first_L_terms} we get that
$$
\COV{F_1(t,h),F_2(t,h)} = 
\frac{9}{\pi \sqrt{\alpha}} 
h^2 \log^2\left( h^{-1}\right).
$$
\end{proof}
\subsubsection{The Case of Rational $\alpha$}\label{sec:alpha_rational}

Assume now that $\alpha\in\mathbb{Q}$ and denote $\alpha = p/q$ where $p,q$ are co-prime.

\begin{lemma}\label{eq:sec4sub2_rational_diophantine}
The pair $(\Omega_1,\Omega_2)$ is a $1$-Diophantine pair.
\end{lemma}

\begin{proof}
Let $n,m\in\ZZ^2$, $0<|n|,|m|<M$.
From \eqref{eq:sec4sub2eq1} we get that if $Y_1(n)\neq Y_2(m)$ then
$$
\left|Y_1(n) - Y_2(m) \right|
\gg 
\frac{1}{M}\left| n_1^2 -  n_1 n_2 +  n_2^2 - m_2^2 - \alpha m_1^2 \right|.
$$
Writing $\alpha = p/q$, we get that
$\left|Y_1(n) - Y_2(m) \right| \gg_q \frac{1}{M}$.
\end{proof}

\begin{proof}[Proof of \autoref{thm:specific_ellipse}, Parts 2,3]
Let $\alpha = \frac{p}{q}\in \mathbb{Q}$. From \autoref{eq:sec4sub2_rational_diophantine} we have that $\Omega_1,\Omega_2$ are a $1$-Diophantine pair.
Thus, from \autoref{thm:covar_short_interval_F1_F2}, if $h(T) \gg T^{-1/62}$ then
\begin{equation*}
\COV{F_1(t,h),F_2(t,h)} 
= 
4\sum_{\nu}
\hat{F_1}(\nu)\overline{\hat{F_2}}(\nu)
\sin^2\left(\frac{h}{2}\nu\right).
\end{equation*}

We denote
$$
r_{2}(\nu) = 
\# \left\{ 
(n_1,n_2)\in\ZZ^2 \quad \middle|  \quad
2\pi \sqrt{ m_1^2 + \frac{q}{p} m_2^2} = \nu
\right\}.
$$

Note that the common frequencies of $F_1$ and $F_2$ are supported on the set $\left\{\pm2\pi \sqrt{\ell} : \ell \in \NN\right\}$, and so we get that
\begin{equation}\label{eq:cov_rational_sec4}
\COV{F_1(t,h),F_2(t,h)} 
= 
\frac{\sqrt{3}}{\pi^2\sqrt{\alpha}}
\sum_{\ell \in \NN}
\frac{r_1\left(2\pi\sqrt{\ell}\right) r_2\left(2\pi\sqrt{\ell}\right)}
{\ell^{3/2}}
\sin^2\left(\pi h \sqrt{\ell}\right).
\end{equation}

In the appendix (Theorems \ref{thm:app_square_case},\ref{thm:app_non_square_case}) we show that
\begin{align}
\sum_{\ell \leq N}r_1\left(2\pi\sqrt{\ell}\right) r_2\left(2\pi\sqrt{\ell}\right)= 
\begin{cases}
    \frac{3\sqrt{3}}{\sqrt{p'q}} N\log N  + \BigO{N} & \text{if } 3pq \text{ is a square }\\
    C  N  +\BigO{N^{3/4 + \epsilon}}             & \text{otherwise}
\end{cases} \label{eq:sub4_cases_rat}
\end{align}
where $p'$ is the squarefree part of $p$, and $C$ is given in \autoref{thm:app_non_square_case}.

We now use \eqref{eq:cov_rational_sec4}, \eqref{eq:sub4_cases_rat} and summation by parts.
We denote $u(x) = \frac{\sin^2\left(\pi \sqrt{x}\right)}{x^{3/2}}$.
We consider the case of $3pq$ not a square.
In this case we get
\begin{multline*}
\COV{F_1(t,h),F_2(t,h)} 
=
\frac{\sqrt{3}}{\pi^2\sqrt{\alpha}}
\sum_{\ell \in \NN}
\frac{r_1\left(2\pi\sqrt{\ell}\right) r_2\left(2\pi\sqrt{\ell}\right)}
{\ell^{3/2}}
\sin^2\left(\pi h \sqrt{\ell}\right)
= \\
\frac{\sqrt{3}}{\pi^2\sqrt{\alpha}}
\sum_{\ell \in \NN}
h^5
\left(C \ell + \BigO{\ell^{3/4+ \epsilon} }\right)
\frac{u\left((\ell+1)h^2\right) - u(\ell h^2)}{h^2}.
\end{multline*}
As $h$ tends to 0, using the definitions of the derivative and of the Riemann sum, we get that
\begin{equation*}
\COV{F_1(t,h),F_2(t,h)} 
=
Ch\frac{\sqrt{3}}{\pi^2\sqrt{\alpha}}
 \int_{0}^{\infty}x u'(x) dx + \BigO{h^{5/4 - \epsilon}}.
\end{equation*}
Using integration by parts and the variable change $y = \sqrt{x}$, it follows that
\begin{multline*}
\COV{F_1(t,h),F_2(t,h)} 
=\\
2Ch\frac{\sqrt{3}}{\pi^2\sqrt{\alpha}}
\int_{0}^{\infty}\frac{\sin^2(\pi y)}{y^2}dy
 = 2Ch\frac{\sqrt{3}}{\pi^2\sqrt{\alpha}} \frac{\pi^2}{2}
 = \frac{C \sqrt{3}}{\sqrt{\alpha}} h .
\end{multline*}

Similar arguments show that in the case where $3pq$ is a square:
\begin{equation*}
\COV{F_1(t,h),F_2(t,h)} = 
\frac{ 18}{\sqrt{p p'}} h \log\left(h^{-1}\right) 
\end{equation*}
where $p'$ is is the squarefree part of $p$ (which is either 3 or 1).

\end{proof}

\appendix
\section{Appendix}\label{sec:appendix}

We denote
\begin{gather*}
 r_\omega(k) = \#\left\{
n_1,n_2\in\ZZ \middle| n_1^2 - n_1n_2 + n_2^2 = k
\right\}   \\[12pt]
r_{a,b}(k) = \# \left\{
n_1,n_2\in\ZZ \middle| n_1^2 + \frac{a}{b}n_2^2 = k 
\right\}.
\end{gather*}
In this section we prove several results regarding these functions that were used in \autoref{thm:specific_ellipse}.
We start by stating the main theorems.

\begin{theorem}\label{thm:mult_case}
$$
\sum_{k\leq N} r_\omega\left(k^2\right) =
\frac{3\sqrt{3}}{\pi}N\log N + \BigO{N}.
$$
\end{theorem}

\begin{theorem}\label{thm:app_square_case}
Assume that $3ab$ is an integer square.
Then
$$
\sum_{k\leq N} 
r_\omega\left(k\right)
r_{a,b}\left(k\right) 
=
\frac{3\sqrt{3}}{\sqrt{ab'}} N\log N
+ \BigO{N}
$$
where $b'$ is the squarefree part of $b$.
\end{theorem}
\begin{theorem}\label{thm:app_non_square_case}
Assume that $3ab$ is not an integer square.
Then
$$
\sum_{k\leq N} 
r_\omega\left(k\right)
r_{a,b}\left(k\right) 
=
C N + \BigO{N^{3/4 + \epsilon}}
$$
for any $\epsilon>0$.

The constant $C$ is given by 
$$
C =\left( \frac{2\pi^2}{\sqrt{3ab'}}\right)
\left(\frac{ L(1,\chi) }{L(2,\chi)}\right)
\sigma_2
\prod_{\substack{p^r \mid\mid 3ab' \\ p\geq3, r\geq 2} } 
\sigma_p 
$$
where-
\begin{itemize}
    \item b' is the squarefree part of b.
    \item  $\chi$ is the Kronecker symbol $\chi(*) = \left(\frac{12ab'}{*}\right)$.
    \item The value of $\sigma_2$ is given in \autoref{lem:a-9}:
    \begin{align*}
\sigma_2 =     
\begin{cases}
    1.5 & \frac{ab'}{4^s} \equiv 3\mod 8\\
    1.5 - \frac{1}{3\cdot 2^s} & \frac{ab'}{4^s} \equiv 7 \mod 8\\
    1.5 - \frac{1}{2^{s+1}} &  \frac{ab'}{4^s} \equiv 1,2,5,6 \mod 8
\end{cases}
\end{align*}
    where $4^s \mid ab'$ and $4^{s+1} \nmid ab'$.
    \item The value of $\sigma_3$ is given in \autoref{lem:a-11}:
    \begin{align*}
\sigma_3 =     
\begin{cases}
    2 - \frac{1}{3^s} & \frac{ab'}{9^s} \not\equiv 0\mod 3\\
    2  & \frac{ab'}{3^{2s+1}} \equiv 1 \mod 3\\
    2 - \frac{1}{2\cdot 3^s} &  \frac{ab'}{3^{2s+1}} \equiv 2 \mod 3
\end{cases}
\end{align*}
    where $9^s \mid ab'$ and $9^{s+1} \nmid ab'$.
    \item The value of $\sigma_p$ for $p> 3$ which divides $ab'$ with even multiplicity $r$ is given in lemma \ref{lem:a-6}:
    \begin{align*}
    \sigma_p = 
    \begin{cases}
        1 + \frac{1}{p} & \frac{3ab'}{p^r} \text{ is a square mod } p\\
        1 + \frac{1}{p} - \frac{2}{p^{r/2}(p+1)} & \text{otherwise}.
    \end{cases}
\end{align*}
    \item The value of $\sigma_p$ for $p> 3$ which divides $ab'$ with odd multiplicity $r$ is given in lemma \ref{lem:a-7}:
$$
\sigma_p = 1 + \frac{1}{p} - \frac{1}{p^{(r+1)/2}}.
$$
\end{itemize}
\end{theorem}

We begin with a proof of \autoref{thm:mult_case}.
\begin{proof}
Let $\omega = e^{2\pi i / 3}$ and let $\ZZ[\omega]$ be the Eisenstein integers.
Denote also $U = \left\{ \pm 1, \pm \omega, \pm \omega^2\right\}$ the units of $\ZZ[\omega]$.
Then $r_\omega(k^2)$ is the number of elements in $\ZZ[\omega]$ with norm $k^2$.
Denote $a(k) = \frac{1}{6}r_\omega(k^2)$.
Then $a(k)$ counts the number of Eisenstein integers of norm $k^2$ up to multiplication by a unit in $U$.
From unique factorization in $\ZZ[\omega]$ we have that $a(k)$ is a multiplicative function.
From properties of primes in $\ZZ[\omega]$ we get that $a\left(p^e\right)$ is:

\begin{align}
a\left(p^e\right) =
\begin{cases}
    2e+1 &   p \equiv 1 \mod 3   \\
    1   &   p \equiv 2 \mod 3   \\
    1   &   p = 3   .
\end{cases} 
\end{align}

By the Euler product expansion, it can be seen that the Dirichlet series of $a\left(k\right)$ is given by
$$
D_a(s) = \sum_{n\in \NN}\frac{a(n)}{n^s}
=
\zeta^2(s) \frac{L(s,\chi_3)}{\zeta(2s)\left(1-3^{-s}\right)}.
$$
We can now use \cite[Theorem 1.1]{S-D}.
This Theorem gives an asymptotic expression for the sum $\sum_{n\leq N} a(n)$, where the corresponding Dirichlet series $D_a(s) = \sum_{n\in \NN}\frac{a(n)}{n^s}$ satisfies $D_a(s) = G(s) \zeta^z(s) \zeta^{-w}(2s)$
where $G(s)$ can be continued as a holomorphic function to some suitable neighborhood of $\Re(s)\geq 1$ with suitable bounds. 
Applied to our case, we get that
$$
\sum_{n\leq N }a\left(k\right) = C N\log N + \BigO{N}
$$
with
$$
C = \frac{L(1,\chi_3)}{\zeta(2)(1 - \frac{1}{3})} = \frac{\sqrt{3}}{2\pi}.
$$
The result follows after noting that $r_\omega(k) = 6 a(k)$.
\end{proof}

We now present the proofs of \autoref{thm:app_square_case} and \autoref{thm:app_non_square_case}.
These use results regarding the singular integral and singular series, which we define and prove in lemmas \ref{lem:singular_integral} through \ref{lem:a-11}.
The terms singular integral and singular series are borrowed from \cite{h-b-circle-method}.

\begin{proof}[Proof of \autoref{thm:app_square_case}]
The sum
$$
S_N =
\sum_{k\leq N} 
r_\omega\left(k\right)
r_{a,b}\left(k\right)
$$
is equal to the number of solutions $x,y,z,t\in\ZZ^4$ of $G(x,y,z,t)= 0$ with the quadratic form
$$
G(x,y,z,t) = x^2 - xy + y^2 -z^2 -\frac{a}{b} t^2.
$$
which satisfy $x^2 - xy + y^2 \leq N$.
Denote by $b'$ the squarefree part of $b$.
Then any solution to $G(x,y,z,t) = 0$ must have $\left(b'\sqrt{\frac{b}{b'}}\right)\mid t$.
Writing $t =\left(b'\sqrt{\frac{b}{b'}}\right) w$ we get that counting solutions to $G=0$ is equivalent to counting solutions $x,y,z,w\in\ZZ^4$ to $F(x,y,z,w)= 0$ where $F$ is the quadratic form
$$
F(x,y,z,w) = x^2 - xy + y^2 -z^2 -a b' w^2.
$$

The discriminant of the quadratic form $F$ is $12 a b'$.
If we assume that $3a b'$ is a square, then we can use \cite[Theorem 7]{h-b-circle-method}.
This theorem gives a formula for the asymptotics of 
$$
\mathcal{N}(P) = \#\left\lbrace F(x,y,z,w) = 0 \mid (x,y,z,w) \in P \mathcal{R}\right\rbrace
$$
where $\mathcal{R}$ is a some compact region, and $F$ is a quadratic form with square discriminant.
This can be applied to our sum $S_N$, since 
$$
S_N = 
\sum_{k\leq N} 
r_\omega\left(k\right)
r_{a,b}\left(k\right)
= 
\#\left\lbrace F(x,y,z,w) = 0 \mid (x,y,z,w)\in N \mathcal{R}\right\rbrace
$$
where
$$
\mathcal{R} = \left\lbrace x^2 -xy + y^2 \leq 1, z^2 + ab'w^2 \leq 1\right\rbrace.
$$
Applying \cite[Theorem 7]{h-b-circle-method} to $S_N$ we get that 
$$
S_N = C N\log N + \BigO{N}.
$$
Furthermore, the constant $C$ is given by
\begin{equation}\label{eq:sing_series_int_for_const_square}
C = \frac{1}{2}\sigma_\infty 
\prod_{p}\left(\left(1-\frac{1}{p}\right) \sigma_p\right)
\end{equation}
where $\sigma_\infty$ is the singular integral, and $\sigma_p$ are the $p$-adic densities.
Since $3\alpha$ is a square, from lemmas \ref{lem:a-5}, \ref{lem:a-6}, \ref{lem:a-9}, \ref{lem:a-11}, we get that $\sigma_p = 1+\frac{1}{p}$ for all $p\neq 3$, and
$\sigma_3 = 2$.
Also, from \autoref{lem:singular_integral} we get that $\sigma_\infty = \frac{2\pi^2}{\sqrt{3ab'}}$.
Plugging these into \eqref{eq:sing_series_int_for_const_square}, we get that
$$
C = \frac{1}{2}
\left(\frac{2\pi^2}{\sqrt{3ab'}}\right)
\left(\frac{6}{\pi^2}\right)
\left(\frac{2}{4/3}\right) =
\frac{3\sqrt{3}}{\sqrt{ab'}}.
$$
\end{proof}

\begin{proof}[Proof of \autoref{thm:app_non_square_case}]
Once more, the sum
$$
S_N =
\sum_{k\leq N} 
r_\omega\left(k\right)
r_{a,b}\left(k\right)
$$
is equivalent to the  number of solutions $x,y,z,w\in\ZZ^4, F(x,y,z,w)= 0$ to the quadratic form
$$
F(x,y,z,w) = x^2 - xy + y^2 -z^2 -a b' w^2
$$
where $b'$ is the squarefree part of $b$.
If we assume that $3ab'$ is not a square then we can use \cite[Theorem 6]{h-b-circle-method}. 
This theorem is similar to \cite[Theorem 7]{h-b-circle-method} except it deals with the case where the discriminant of the quadratic form $F$ is not a square.
Applied to $S_N$ we get that
$S_N = C N + \BigO{N^{3/4 + \epsilon}}$ for any $\epsilon > 0$.
Furthermore, we get a formula for the constant $C$.
Denote by $\chi$ the Kronecker symbol
$\chi(*) = \left(\frac{disc(F)}{*}\right)$ where $disc(F) = 12ab'$.
Then
\begin{equation*}
C = \sigma_\infty L(1,\chi) 
\prod_{p}\left(1 - \chi(p) p^{-1}\right)\sigma_p.
\end{equation*}
From lemmas \ref{lem:a-5}, \ref{lem:a-7} \ref{lem:a-11} we get that
\begin{align*}
    \sigma_p = 
    \begin{cases}
        1 + \frac{1}{p} & \chi(p) = 1 \\
        \left(1 + \frac{1}{p^2}\right)
        \left(1 + \frac{1}{p}\right) & \chi(p) = -1 \\
        1 &     p\neq 2 \quad p \mid 3a'b \text{ and } p^2 \nmid 3a'b.
    \end{cases}
\end{align*}
The cases of $p=2$ and $p^r \mid 3ab'$ for some $r>1$ are handled in lemmas \ref{lem:a-6},\ref{lem:a-7},\ref{lem:a-9},\ref{lem:a-11}.
It follows that
\begin{equation*}
\prod_{p}\left(1 - \chi(p) p^{-1}\right)\sigma_p=
 L(2,\chi)^{-1} 
\prod_{\substack{p^r \mid\mid 12ab' \\  r\geq 2} } 
\sigma_p.
\end{equation*}

where the values of $\sigma_p$ for $p=2$ and $p\geq3$ with $p^r\mid\mid 3ab'$ for some $r\geq 2$ are explicitly given in lemmas \ref{lem:a-6},\ref{lem:a-7},\ref{lem:a-9},\ref{lem:a-11}.

As for the singular integral, from \autoref{lem:singular_integral} we get $\sigma_\infty = \frac{2\pi^2}{\sqrt{3ab'}}$.
Thus, overall we get that
$$
C =\left( \frac{2\pi^2}{\sqrt{3ab'}}\right)
\left(\frac{ L(1,\chi) }{L(2,\chi)}\right)
\sigma_2
\prod_{\substack{p^r \mid\mid 3ab' \\ p\geq3, r\geq 2} } 
\sigma_p .
$$
For example, if $3ab'$ is squarefree and $ab'\equiv 1,5 \mod 8$, then we get that
$$
C =\left( \frac{2\pi^2}{\sqrt{3a'b}}\right)
\left(\frac{ L(1,\chi) }{L(2,\chi)}\right).
$$
\end{proof}

We now turn to the calculations of the singular integral and $p$-adic densities.
Throughout, we denote
$$
F(x,y,z,w) = x^2 - xy + y^2 - z^2 -ab'w^2
$$
where $\gcd(a,b')=1$ and $b'$ is squarefree.
We will at times denote 
$$\alpha = ab' .$$

\begin{lemma}\label{lem:singular_integral}
Let 
\begin{align*}
R(x,y,z,w) = 
\begin{cases}
    1&   x^2 -xy + y^2  \leq 1   \\
    0&   \text{ otherwise}.
\end{cases} 
\end{align*}
Then the singular integral $\sigma_\infty$, which is given by
$$
\sigma_\infty = 
\lim_{\epsilon\rightarrow 0^+}
\frac{1}{2\epsilon}
\int_{\left|F(x,y,z,w)\right|< \epsilon}
R(x,y,z,w)dxdydzdw,
$$
has value
$$
\sigma_\infty = \frac{2\pi^2}{\sqrt{3\alpha}}.
$$
\end{lemma}

\begin{proof}
We first consider the following coordinate change:
$$
(x,y,z,w) \mapsto 
\left(
x',y',z',w'
\right) = 
\left(
\frac{x-2y}{2} , 
\frac{\sqrt{3}}{2} x, 
 z, 
\sqrt{\alpha} w
\right)
$$
with Jacobian
$$
\frac{\partial(x',y',z',w')}{\partial(x,y,z,w)} = \frac{\sqrt{3\alpha}}{2}.
$$
Note that
$$
F(x',y',z',w') = x'^2 + y'^2 - z'^2 -w'^2.
$$
We then have 
$$
\sigma_\infty = 
\frac{2}{\sqrt{3\alpha}}
\lim_{\epsilon\rightarrow 0}
\frac{1}{2\epsilon}
\text{Vol}\left(V_\epsilon\right)
$$
where 
$$
V_\epsilon = \left\{
(x',y',z',w') \middle|
\left|{x'}^2 + {y'}^2 - {z'}^2 -{w'}^2\right|\leq \epsilon, \quad
{x'}^2 + {y'}^2 \leq 1
\right\}.
$$

We now switch to polar coordinates in both ${x'}{y'}$-plane and ${z'}{w'}$-plane:
\begin{gather*}
    (x',y')\mapsto (r,\theta),\\
    (z',w')\mapsto (\rho,\psi).
\end{gather*}

We then get
\begin{multline*}
\frac{1}{2\epsilon}
\text{Vol}\left(V_\epsilon\right)=
\frac{1}{\epsilon}
\int_{\theta = 0}^{2\pi} d\theta 
\int_{\psi = 0}^{2\pi} d\psi
\int_{r=0}^{\sqrt{1-\epsilon}}r
\int_{\rho = r}^{\sqrt{r^2 + \epsilon}} \rho d\rho dr
=\\
\frac{4\pi^2}{\epsilon}
\left(\frac{\rho^2}{2}\right)\Big|_{r}^{\sqrt{r^2 + \epsilon}}
\left(\frac{r^2}{2}\right)\Big|_{0}^{\sqrt{1 - \epsilon}}
\xrightarrow[]{\epsilon\rightarrow 0} \pi^2.
\end{multline*}
It follows that 
$$
\sigma_\infty = \frac{2\pi^2}{\sqrt{3\alpha}}.
$$
\end{proof}

We now proceed to compute the $p$-adic densities:
$$
\sigma_p = \lim_{k\rightarrow\infty}\frac{1}{p^{3k}}
\#\left\{
(x,y,z,w)\quad  (\text{mod} p^k) \middle| F(x,y,z,w) = 0 \quad(\text{mod} p^k)
\right\}.
$$
We first introduce some notations.
Let $p$ be a prime.
Assume that $p^i \mid \alpha$ and denote $\alpha' = \frac{ab'}{p^i}$.
We denote
$$
N_k^r = \#\left\{
(x,y,z,w)\text{ }  (\text{mod} p^k) \middle| 
x^2 - xy + y^2 - z^2 -p^r\alpha' w^2 = 0
\text{ }(\text{mod} p^k)
\right\}.
$$
We also denote $N_k = N_k^0$, and we define $N_0 = 1$.

\begin{lemma}\label{lem:a-5}
Let $p$ be a prime, $p\neq 2,3$, $p\nmid \alpha$ (where $\alpha = ab'$).
Then 
\begin{align*}
    \sigma_p = 
    \begin{cases}
        1 + \frac{1}{p} &  3\alpha \text{ is a square in } \mathbb{F}_p  \\
        \left(1 + \frac{1}{p^2}\right)\left(1 + \frac{1}{p}\right)^{-1} & \text{otherwise}.
    \end{cases}
\end{align*}

\end{lemma}
\begin{proof}
Since $p\neq2,3$, and
$$
x^2 -xy + y^2 = \frac{1}{4}\left((x-2y)^2 + 3x^2\right),
$$
we can consider the equivalent form $x^2 + 3y^2 - z^2 - \alpha w^2$.

We first find the number of solutions mod $p$.
We can consider an equivalent form
$$
F'(x,y,z,h):
(x+h)^2 + 3y^2 - h^2 -\alpha w^2
$$
which can be rewritten as
$$
2xh + x^2 +3y^2 - \alpha w^2.
$$
For any value of $x,y,w\mod p $ for which $x \neq 0$ we can pick $h$ accordingly to solve the equation.
This gives us $(p-1)p^2$ solutions.

If however we assume $x = 0$ mod $p$, then (after picking arbitrary $h$ which has $p$ options) our equation becomes
\begin{equation}\label{eq:aaa_non_linear}
3y^2 = \alpha w^2 \mod p    .
\end{equation}
If $3\alpha$ is a square mod $p$, then this has $2p - 1$ solutions.
Otherwise, there is only the single 0 solution.
And so, we get an additional
$p(2p-1)$ solutions if  $3\alpha $ is a square mod $p$, and only an additional $p$ solutions otherwise.

Thus, the number of non-zero solutions mod $p$ satisfies
\begin{align*}
    N_1 - 1 = 
    \begin{cases}
        p^3 +p^2 -p -1 & 3\alpha \text{ is a square }\mod p\\
        p^3 -p^2 +p - 1 & \text{otherwise}.
    \end{cases}
\end{align*}

We now use Hensel's lemma to calculate the number of solutions mod $p^k$.
We say that a solution mod $p^k$ is non-degenerate if it is not the zero solution mod $p$.
From Hensel's lemma, every non-degenerate solution can be lifted from $p^i$ to $p^{i+1}$ in $p^3$ ways.
This means that we have $p^{3(k-1)}\left(N_1-1\right)$ non-degenerate solutions mod $p^k$.

For degenerate solution mod $p^k$, we can divide each variable by $p$.
This gives us a solution mod $p^{k-2}$, with the most significant p-adic digit of each variable chosen arbitrarily.
It follows that the number of degenerate solutions mod $p^k$ is $p^4 N_{k-2}$.

Thus we have the recursion
\begin{equation*}
\frac{N_k}{p^{3k}}
=
\frac{N_1 - 1}{p^3} + \frac{1}{p^2}\left(\frac{N_{k-2}}{p^{3(k-2)}}\right).
\end{equation*}
This recursion can be solved to give
\begin{equation*}
\lim_{k\rightarrow\infty}\frac{N_k}{p^{3k}} = 
\left(1 +\frac{1}{p} -\frac{1}{p^2} -\frac{1}{p^3}\right)
\left(1 - \frac{1}{p^2}\right)^{-1} 
= 1 + \frac{1}{p}
\end{equation*}
if $3\alpha$ is a square mod $p$, and
\begin{equation*}
\lim_{k\rightarrow\infty}\frac{N_k}{p^{3k}} = 
\left(1 -\frac{1}{p} +\frac{1}{p^2} -\frac{1}{p^3}\right)
\left(1 - \frac{1}{p^2}\right)^{-1} 
=
\left(1 + \frac{1}{p^2}\right)\left(1 + \frac{1}{p}\right)^{-1}
\end{equation*}
otherwise.
\end{proof}

\begin{lemma}\label{lem:a-6}
Let $p\neq 2,3$.
Assume that $p^r\mid \alpha$, where $r>0$ is even and $p^{r+1}\nmid \alpha$ (where $\alpha = ab'$).
Then
\begin{align*}
    \sigma_p = 
    \begin{cases}
        1 + \frac{1}{p} & \frac{3\alpha}{p^r} \text{ is a square mod } p\\
        1 + \frac{1}{p} - \frac{2}{p^{r/2}(p+1)} & \text{otherwise}.
    \end{cases}
\end{align*}
\end{lemma}
\begin{proof}
We consider the equation
\begin{equation}\label{eq:b:0}
x^2 - xy + y^2 - z^2 -p^r \alpha w^2 = 0 \mod p^k.    
\end{equation}

Modulo $p$, this equation become:
\begin{equation}\label{eq:b:1}
x^2 - xy + y^2 - z^2 = 0\mod p
\end{equation}
(and $w$ can be chosen arbitrarily).
We will call a solution mod $p$ good if either $x,y$ or $z$ is inevitable. 
From Hensel's lemma, these good solutions can be lifted to a solution of \eqref{eq:b:0} mod $p^k$ in $p^{3(k-1)}$ ways.
It can be shown that the number of good solutions to \eqref{eq:b:1} is $p^2 - 1$ (times $p$ options for $w$).
So the number of solutions mod $p^k$ which come from a good solution is 
$(p^3 - p)p^{3(k-1)}$.

If $(x_0,y_0,z_0,w_0)$ is a solution that doesn't come from a good solution, then we can divide $x_0,y_0,z_0$ by $p$ and the coefficient of $w_0$ by $p^2$ and get a solution mod $p^{k-2}$.
We can then pick the most significant p-adic digit of $x_0,y_0,z_0$ arbitraraly, and the two most significant p-adic digit of $w_0$.
Thus, the number of non-good solutions is $p^5 N_{k-2}^{r-2}$.

Overall, we get the following recursion:
\begin{equation}\label{eq:b:2}
\frac{N_k^r}{p^{3k}} = 
\left(1 - \frac{1}{p^2}\right)
+ \frac{1}{p}\left(\frac{N_{k-2}^{r-2}}{p^{3(k-2)}}\right).
\end{equation}

By repeating \eqref{eq:b:2} $r/2$ times we get
\begin{multline}\label{eq:b:3}
\frac{N_k^r}{p^{3k}} = 
\left(1-\frac{1}{p^2}\right)
\left( 1+ \frac{1}{p} + \frac{1}{p^2} + ... + \frac{1}{p^{r/2 - 1}}\right)
+ \frac{1}{p^{r/2}}\left(\frac{N_{k-r}}{p^{3(k-r)}}\right)
=\\
\left(1-\frac{1}{p^2}\right)
\frac{1 - \frac{1}{p^{r/2}}}{1 - \frac{1}{p}}
+\frac{1}{p^{r/2}} \left( \frac{N_{k-r}}{p^{3(k-r)}}\right)
= \\
\left(1 + \frac{1}{p}\right)
+\frac{1}{p^{r/2}}\left(\frac{N_{k-r}}{p^{3(k-r)}} - 1 - \frac{1}{p}\right).
\end{multline}

From \autoref{lem:a-5}, we have
\begin{align*}
\lim_{k\rightarrow\infty} \frac{N_{k-r}}{p^{3(k-r)}} 
= 
\begin{cases}
    1 + \frac{1}{p} & 3\alpha \text{ is a square mod } p\\
     \left(1 + \frac{1}{p^2}\right)\left(1 + \frac{1}{p}\right)^{-1} & \text{otherwise}.
\end{cases}
\end{align*}

Plugging this into \eqref{eq:b:3} we get 

\begin{align*}
\lim_{k\rightarrow\infty} \frac{N_k^r}{p^{3k}} = 
= 
\begin{cases}
    1 + \frac{1}{p} & 3\alpha \text{ is a square mod } p\\
    1 + \frac{1}{p} - \frac{2}{p^{r/2}(p+1)} & \text{otherwise}.
\end{cases}
\end{align*}
\end{proof}

\begin{lemma}\label{lem:a-7}
Let $p\neq 2,3$.
Assume that $p^r\mid \alpha$, $p^{r+1}\nmid \alpha$ where $r$ is odd (where $\alpha = ab'$).
Then
$$
\sigma_p = 1 + \frac{1}{p} - \frac{1}{p^{(r+1)/2}}.
$$
\end{lemma}
\begin{proof}
We start with the case $r=1$.
That is, we show that $\frac{N_k^1}{p^{3k}} \rightarrow 1$.

For this, we first fine the non-degenerate solutions mod $p$.
Our equation mod $p$ is
$$
x^2 + 3y^2 - z^2 = 0 \mod p.
$$
There are then $p(p^2-1)$ solutions $x_0,y_0,z_0,w_0$ with at least one of $x_0,y_0$, or $z_0$ not zero.
From Hensel's lemma, these give us $p(p^2 - 1) p^{3(k-1)}$ solution mod $p^k$.

We now consider the rest of the solutions.
Let $x_0,y_0,z_0,w_0$ be a solution mod $p^k$ such that $x_0,y_0$ and $z_0$ are 0 mod $p$.
Assume also that $k>1$.
Then since
$$
p \alpha w_0^2 = x_0^2 + 3y_0^2 - z_0^2 \mod p^k
$$
We must have that $p\mid w_0$ as well (since the RHS is divisible by $p^2$).
Dividing $x_0,y_0,z_0,w_0$ by $p$ will then give us a solution mod $p^{k-2}$.

This gives us the relation
\begin{equation*}
\frac{N_k^1}{p^{3k}} 
=
\left(1 - \frac{1}{p^2}\right) 
+ 
\frac{1}{p^2}\left(\frac{N_{k-2}^1}{p^{k(3-2)}}\right).
\end{equation*}

From this we can deduce that
\begin{equation*}
\lim_{k\rightarrow\infty}
\frac{N_k^1}{p^{3k}}  = 
\left(1 - \frac{1}{p^2}\right) 
\left( \frac{1}{1 - \frac{1}{p^2}}\right)
= 1.
\end{equation*}

We now consider the case $r>1$.
It can be shown in the same way we did in \autoref{lem:a-6} that
\begin{equation*}
\frac{N_k^r}{p^{3k}} = 
\left(1 - \frac{1}{p^2}\right)
+ \frac{1}{p}\left(\frac{N_{k-2}^{r-2}}{p^{3(k-2)}}\right).
\end{equation*}

Repeating this recursion $(r-1)/2$ times gives us
\begin{multline*}
\frac{N_k^r}{p^{3k}} = 
\left(1-\frac{1}{p^2}\right)
\left( 1+ \frac{1}{p} + \frac{1}{p^2} + ... + \frac{1}{p^{(r-3)/2 - 1}}\right)
+ \frac{1}{p^{(r-1)/2}}\left(\frac{N_{k-r + 1}^1}{p^{3(k-r + 1)}}\right)
\\=
\left(1-\frac{1}{p^2}\right)
\frac{1 - \frac{1}{p^{(r-1)/2}}}{1 - \frac{1}{p}}
+\frac{1}{p^{(r-1)/2}} \left( \frac{N_{k-r+1}^1}{p^{3(k-r+1)}}\right)
\\= 
\left(1 + \frac{1}{p}\right)
+\frac{1}{p^{(r-1)/2}}\left(\frac{N_{k-r+1}^1}{p^{3(k-r+1)}} - 1 - \frac{1}{p}\right).
\end{multline*}

Using the fact that 
$$
\frac{N_{k-r+1}^1}{p^{3(k-r+1)}} \xrightarrow[]{k\rightarrow\infty} 1
$$ 
gives us
\begin{equation*}
    \frac{N_k^r}{p^{3k}} \xrightarrow[]{k\rightarrow\infty}
    1 + \frac{1}{p} - \frac{1}{p^{(r+1)/2}}.
\end{equation*}
\end{proof}

\begin{lemma}\label{lem:app_sigular_2_ez_case}
Assume that $4\nmid \alpha$ (where $\alpha = ab'$), then
\begin{align*}
\sigma_2 =     
\begin{cases}
    1.5 & \alpha \equiv 3\mod 8\\
    7/6 & \alpha \equiv 7 \mod 8\\
    1 &  \alpha \equiv 1,2,5,6 \mod 8.
\end{cases}
\end{align*}
\end{lemma}

\begin{proof}
The difference from previous case is that we can only use Hensel's lemma starting at $2^3$.
Thus, every non-degenerate solution mod $2^3$ can be lifted to a non-degenerate solution mod $2^k$ in $2^{3(k-3)}$ ways.
This is also true if $2\mid \alpha$, $4\nmid \alpha$ since then every non-degenerate solution mod $2^3$ has to have one of $x,y$ or $z$ be inevitable.

Denote by $n_3$ the number of non-degenerate solutions mod $2^3$.
In a similar way to previous lemmas, we get the recursion-
\begin{equation*}
\frac{N_k}{2^{3k}} = 
\frac{n_3}{8^3} + 
\frac{1}{4}
\left(\frac{N_{k-2}}{2^{3(k-2)}}\right).
\end{equation*}
This implies
$$
\frac{N_k}{2^{3k}} \xrightarrow[]{k\rightarrow\infty} \frac{4}{3}\left(\frac{n_3}{8^3}\right).
$$
It remains to find the value of $\frac{n_3}{8^3}$.
This can be done by direct computation:
\begin{align*}
\frac{n_3}{8^3} = 
\begin{cases}
    1.125 & \alpha \equiv 3\mod 8\\
    0.875 & \alpha \equiv 7 \mod 8\\
    0.75 &  \alpha \equiv 1,2,5,6 \mod 8.
\end{cases}    
\end{align*}
The result then follows.
\end{proof}

\begin{lemma}\label{lem:a-9}
If $4^s \mid \alpha$, $4^{s+1}\nmid \alpha$ (where $\alpha = ab'$), then
$$
\sigma_2 = 1.5 + \frac{1}{2^{s}}\left(\sigma_2 ' - 1.5\right)
$$
where $\sigma_2 '$ is the value of $\sigma_2$ from \autoref{lem:app_sigular_2_ez_case} for $\frac{\alpha}{2^{2s}}\mod 8$.
\end{lemma}
\begin{remark}
If $3\alpha$ is a square in the integers, then $\frac{\alpha}{2^{2s}}\equiv 3 \mod 8$.
This means that $\sigma_2 ' = 1.5$.
Thus, for this case we get $\sigma_2 = 1.5$.
\end{remark}

\begin{proof}
In similar ways to previous lemmas, we get the recursion
\begin{equation}\label{eq:rec_case2}
\frac{N_k^r}{2^k}
=
\frac{n_3}{8^3} + 
\frac{1}{2}
\left(\frac{N_{k-2}^{r-2}}{2^{3(k-2)}}\right).
\end{equation}
The result then follows in a similar way to \autoref{lem:a-7},
by repeating \eqref{eq:rec_case2} $s = \lfloor r/2 \rfloor$ times and then using \autoref{lem:app_sigular_2_ez_case}. 
\end{proof}

\begin{lemma}\label{lem:a-10}
Assume that $3^2 \nmid \alpha$ (where $\alpha = ab'$).
Then-
\begin{align*}
\sigma_3 = 
\begin{cases}
    1       &   3 \nmid \alpha \\
    2       &   \frac{\alpha}{3} \equiv 1 \mod 3 \\
    1.5     &   \frac{\alpha}{3} \equiv 2 \mod 3 .
\end{cases}
\end{align*}
\end{lemma}

\begin{proof}
Mod $3^k$, the quadratic form $F$ is equivalent to
$$
x^2 + 3 y^2 - z^2 - \alpha w^2.
$$

Assuming that $3^2\nmid \alpha$, we can use Hensel's starting at $3^3$.
Denote by $n$ the number of non-degenerate solutions mod $3^3$.
We can get the recursion
$$
\frac{N_k}{3^{3k}} = \frac{n}{3^3} + 
\frac{1}{p^2}\left( \frac{N_{k-2}}{3^{3(k-2)}}\right).
$$
From this we get that
$$
\lim_{k\rightarrow\infty }\frac{N_k}{3^{3k}} = 
\frac{9}{8} \left(\frac{n}{3^3}\right).
$$
From direct computation we get that
\begin{align*}
    \frac{9}{8} \left(\frac{n}{3^3}\right) = 
    \begin{cases}
        1   &   3\nmid \alpha \\
        2   &   3\mid \alpha, \quad \frac{\alpha}{3} \equiv 1 \mod 3 \\
        1.5 &    3\mid \alpha, \quad\frac{\alpha}{3} \equiv 2 \mod 3 
    \end{cases}
\end{align*}
\end{proof}

\begin{lemma}\label{lem:a-11}
Assume $9^s \mid \alpha$, $9^{s+1} \nmid \alpha$ for some $s\geq0$ (where $\alpha = ab'$), then
$$
\sigma_3 = 2 - \frac{1}{3^s}\left(2 - \sigma_3'\right)
$$
where $\sigma_3'$ is the value from \autoref{lem:a-10} corresponding to $\frac{\alpha}{9^s}$.
\end{lemma}

\begin{proof}
We have the equation 
\begin{equation}\label{eq:a10-1}
x^2 + 3y^2 - z^2 - \alpha w^2 = 0 \mod 3^k    
\end{equation}
for some $k\geq 3$.
We will call a solution $x_0,y_0,z_0,w_0$ of \eqref{eq:a10-1} good if at least one of $x_0,y_0$ or $z_0$ is invertible.
From Hensel's lemma, every good solution mod $3^3$ can be lifted in $3^{3(k-3)}$ ways to a solution mod $p^k$.
From direct calculation we get that the number $n$ of good solutions mod $3^3$ satisfies $\frac{n}{3^3} = \frac{4}{3}$.

Every non-good solution comes from a solution in $N_{k-2}^{r-2}$ in $p^5$ ways.
This gives us the recursion
$$
\frac{N_k^r}{3^k} = \frac{4}{3} + 
\frac{1}{3}\left(\frac{N_{k-2}^{r-2}}{3^{3(k-2)}}\right).
$$

Repeating this recursion $s = \lfloor r/2 \rfloor$ times and using \autoref{lem:a-10} we get
$$
\sigma_3 = 2 - \frac{1}{3^s}\left(2 - \sigma_3'\right)
$$
where $\sigma_3'$ is the value from \autoref{lem:a-10} corresponding to $\frac{\alpha}{9^s}$.

\end{proof}


\begin{thebibliography}{10}

\bibitem{berry-tabor}
Michael~Victor Berry and Michael Tabor.
\newblock Level clustering in the regular spectrum.
\newblock {\em Proceedings of the Royal Society of London. A. Mathematical and
  Physical Sciences}, 356(1686):375--394, 1977.

\bibitem{besicovitch}
Abram~Samoilovitch Besicovitch.
\newblock {\em Almost periodic functions}, volume~4.
\newblock Dover New York, 1954.

\bibitem{bleher}
Pavel Bleher et~al.
\newblock On the distribution of the number of lattice points inside a family
  of convex ovals.
\newblock {\em Duke Mathematical Journal}, 67(3):461--481, 1992.

\bibitem{bleher_lebowitz}
Pavel~M. Bleher and Joel~L. Lebowitz.
\newblock Variance of number of lattice points in random narrow elliptic strip.
\newblock {\em Annales de l'I.H.P. Probabilit\'es et statistiques},
  31(1):27--58, 1995.

\bibitem{B-W}
Jean Bourgain and Nigel Watt.
\newblock Mean square of zeta function, circle problem and divisor problem
  revisited.
\newblock {\em arXiv preprint arXiv:1709.04340}, 2017.

\bibitem{S-D}
Z Cui, J Wu. 
\newblock The Selberg–Delange method in short intervals with an application. 
\newblock {\em Acta Arithmetica}, 163(3):247--260, 2014.

\bibitem{G-W-W}
Carolyn Gordon, David Webb, and Scott Wolpert.
\newblock Isospectral plane domains and surfaces via riemannian orbifolds.
\newblock {\em Inventiones mathematicae}, 110(1):1--22, 1992.

\bibitem{hardy}
Godfrey~Harold Hardy.
\newblock On the expression of a number as the sum of two squares.
\newblock {\em Quart. J. Math.}, 46:263--283, 1915.

\bibitem{heath}
DR~Heath-Brown.
\newblock The distribution and moments of the error term in the dirichlet
  divisor problem.
\newblock {\em Acta Arithmetica}, 60(4):389--415, 1992.

\bibitem{h-b-circle-method}
DR~Heath-Brown.
\newblock A new form of the circle method, and its application to quadratic
  forms.
\newblock {\em Journal Fur Die Reine Und Angewandte Mathematik - J REINE ANGEW
  MATH}, 1996:149--206, 01 1996.

\bibitem{ivrii}
V~Ya Ivrii.
\newblock Second term of the spectral asymptotic expansion of the
  laplace-beltrami operator on manifolds with boundary.
\newblock {\em Functional Analysis and Its Applications}, 14(2):98--106, 1980.

\bibitem{kac}
Mark Kac.
\newblock Can one hear the shape of a drum?
\newblock {\em The american mathematical monthly}, 73(4P2):1--23, 1966.

\bibitem{lame}
Gabriel Lam{\'e}.
\newblock M{\'e}moire sur la propagation de la chaleur dans les poly{\`e}dres,
  et principalement dans le prisme triangulaire r{\'e}gulier, 1833.

\bibitem{landau}
E~Landau.
\newblock Ueber die gitterpunkte in einem kreise.(ii, mitteilung.).
\newblock {\em Nachrichten von der Gesellschaft der Wissenschaften zu
  G{\"o}ttingen, Mathematisch-Physikalische Klasse}, 1915:161--171, 1915.

\bibitem{McCartin}
Brian~J McCartin.
\newblock {\em Laplacian eigenstructure of the equilateral triangle}.
\newblock Hikari, Limited, 2011.

\bibitem{zeev-QC}
Ze’ev Rudnick.
\newblock What is quantum chaos?
\newblock {\em Notices of the AMS}, 55(1):32--34, 2008.

\bibitem{sierpinski}
W~Sierpi{\'n}ski.
\newblock O pewnem zagadneniu w rachunku funkcyj asymptoticznych.
\newblock {\em Prace Mat.-Fiz}, 17:77--118, 1906.


\bibitem{weyl}
Hermann Weyl.
\newblock {\"U}ber die asymptotische verteilung der eigenwerte.
\newblock {\em Nachrichten von der Gesellschaft der Wissenschaften zu
  G{\"o}ttingen, Mathematisch-Physikalische Klasse}, 1911:110--117, 1911.

\bibitem{wintner}
Aurel Wintner.
\newblock On the lattice problem of gauss.
\newblock {\em American Journal of Mathematics}, 63(3):619--627, 1941.

\end{thebibliography}
\end{document}